\documentclass[11pt]{amsart} \usepackage{amscd, amssymb}
\usepackage{xypic}
\usepackage[matrix, arrow]{xy}
\usepackage{array}

\begin{document}

\def\A{\mathcal{A}}
\def\H{\mathcal{H}}
\def\h{\mathcal{H}}
\def\E{\mathcal{E}}
\def\G{\mathcal{G}}
\def\F{\mathcal F}
\def\X{\mathcal X}
\def\mor{\operatorname{mor}}
\def\spin{\operatorname{spin}}
\def\hot{\hat{\otimes}}

\newcommand{\Cc}{\mathcal{C}}

\newcommand*{\Comp}{{\mathbb{K}}}

\newcommand{\T}{\mathbb{T}}
\newcommand{\Q}{\mathbb{Q}}
\newcommand{\Z}{\mathbb{Z}}
\newcommand{\C}{\mathbb{C}}
\newcommand{\N}{\mathbb{N}}
\newcommand{\R}{\mathbb{R}}
\newcommand{\Ktimes}{\times\!\!\!\!\!\times}

\newcommand{\Lambdah}{\widehat{\Lambda}}
\newcommand{\Deltah}{\widehat{\Delta}}
\newcommand{\ip}[1]{\langle #1 \rangle}

\newcommand{\Kn}{\mathcal{K}_{\nu}}
\newcommand{\Ko}{\mathcal{K}_{\omega}}
\newcommand{\Kob}{\mathcal{K}_{\bar{\omega}}}
\newcommand{\ltg}{\ell^2G}
\newcommand{\tG}{\tilde{G}}
\newcommand{\tN}{\tilde{N}}
\newcommand{\Go}{G_{\omega}}
\newcommand{\Co}{\mathbb{C}_{\omega}}
\newcommand{\ob}{\bar{\omega}}
\newcommand{\Gob}{G_{\bar{\omega}}}
\newcommand{\Cob}{\mathbb{C}_{\bar{\omega}}}
\newcommand{\ltgt}{\ltg_{\tau}}
\newcommand{\talpha}{\tilde{\alpha}}
\newcommand{\ttau}{\tilde{\tau}}
\newcommand{\FUpi}{\mathcal{F}_U\hot_{C_0(X),p_i}\Gamma_0(E)}
\newcommand{\FUpzero}{\mathcal{F}_U\hot_{C_0(X),p_0}\Gamma_0(E)}
\newcommand{\FUpone}{\mathcal{F}_U\hot_{C_0(X),p_1}\Gamma_0(E)}

\newcommand{\thmcite}[1]{{\bfseries\upshape \cite{#1}}}
\newcommand{\thmcitemore}[2]{{\bfseries\upshape \cite[#2]{#1}}}
\newcommand{\citemore}[2]{{\upshape \cite[#2]{#1}}}

\def\Cl{\mathcal{C}l}
\def\Aut{\operatorname{Aut}}
\def\Gx{X\rtimes G}
\def\Bott{\operatorname{Bott}}
\def\CC{\mathbb C}
\def\ClV{C_{\!{}_{V}}}
\def\comp{\operatorname{comp}}
\def\ddS{\stackrel{\scriptscriptstyle{o}}{S}}
\def\intW{\stackrel{\scriptscriptstyle{o}}{W}}
\def\intU{\stackrel{\scriptscriptstyle{o}}{U}}
\def\intT{\stackrel{\scriptscriptstyle{o}}{T}} \def\E{\mathcal E}
\def\EE{\mathbb E}\def\Egx{\mathcal{E}(G)\times X}
\def\EGG{\mathcal{E}(\mathcal{G})}
 \def\EGT{\mathcal{E}(G/G_0)}\def
\EG{\underline{EG}}
\def\EGN{\mathcal{E}(G/N)}
\def\EH{\mathcal{E}(H)}
\def\I{{\operatorname{I}}}
\def\i{{\operatorname{i}}}
 \def\Ind{\operatorname{Ind}}
\def\Id{\operatorname{Id}}
\def\infl{\operatorname{inf}}
\def\k{\operatorname{K}}
\def\KK{\operatorname{KK}}
\def\L{\mathcal L} \def\lk{\langle}
\def\rk{\rangle}
\def\NN{\mathbb N}
\def\oplus{\bigoplus}
\def\pt{\operatorname{pt}}
\def\QQ{\mathbb Q}
\def\res{\operatorname{res}}
\def\RKK{\mathcal{R}\!\operatorname{KK}}
\def\RR{\mathbb R}
\def\sm{\backslash} \def\top{\operatorname{top}}
\def\ZM{{\mathcal Z}M}
 \def\ZZ{\mathbb Z} \setcounter{section}{-1}
\def\Inf{\operatorname{Inf}}
\def\Ad{\operatorname{Ad}}
\def\K{\mathcal K}
\def\id{\operatorname{id}}
\def\U{\mathcal U}
\def\PU{\mathcal PU}
\def\BG{\operatorname{BG}}
\def\eps{\epsilon}
\def\ev{\operatorname{ev}}
\def\om{\omega}
\def\Om{\Omega}
\def\F{\mathcal{F}}
\def\TT{\mathbb T}

\def\Flk{\F_{L,K}}
\def\mg{\mu_{G,A}}
\def\mgx{\mu_{\Gx,A}}
\def\ts{\hot}
\def\ga{\gamma}
\def\Br{\operatorname{Br}}
\def\Ab{\operatorname{Ab}}
\def\la{\lambda}
\def\gr{\operatorname{gr}}
\def\br{\Br^G_{\gr}(X)}
\def\mod{\operatorname{mod}}
\def\hato{\hat{\hot}}
\def\hotimes{{\hot}}
\def\hatoX{\hat{\hot}_{C_0(X)}}
\def\Ctau{C_{\tau}(X)}
\def\ctau{C_{\tau}}
\def\Exp{\operatorname{Exp}}

\theoremstyle{plain}
 \newtheorem{thm}{Theorem}[section]
\newtheorem{cor}[thm]{Corollary}
 \newtheorem{lem}[thm]{Lemma}
\newtheorem{prop}[thm]{Proposition} \newtheorem{lemdef}[thm]{Lemma and Definition}

\theoremstyle{definition}
\newtheorem{defn}[thm]{Definition}
\newtheorem{defremark}[thm]{\bf Definition and Remark}
\theoremstyle{remark}
\newtheorem{remark}[thm]{\bf Remark} \newtheorem{ex}[thm]{\bf Example}
\newtheorem{assumption}[thm]{\bf Assumption}
\newtheorem{notations}[thm]{\bf Notation}
\numberwithin{equation}{section} \emergencystretch 25pt
\renewcommand{\theenumi}{\roman{enumi}}
\renewcommand{\labelenumi}{(\theenumi)}
\title[KK-theoretic duality for proper twisted actions]
  {$\KK$-theoretic duality for proper twisted actions}

\author[Echterhoff]{Siegfried
  Echterhoff}
   \address{Westf\"alische Wilhelms-Universit\"at M\"unster,
  Mathematisches Institut, Einsteinstr. 62 D-48149 M\"unster, Germany}
\email{echters@math.uni-muenster.de}

  \author[Emerson]{Heath
  Emerson}
  \address{Department of Mathematics and Statistics,
  University of Victoria,
  PO BOX 3045 STN CSC Victoria,
  B.C.Canada}
  \email{hemerson@math.uvic.ca}

  \author[Kim]{Hyun Jeong
  Kim}
\address{Westf\"alische Wilhelms-Universit\"at M\"unster,
  Mathematisches Institut, Einsteinstr. 62 D-48149 M\"unster, Germany}
\email{hjkim@math.uni-muenster.de}

\begin{abstract} Let $\mathcal{A}$ be a smooth continuous trace algebra, with a Riemannian manifold spectrum $X$,
equipped with a smooth action by a discrete group $G$ such that $G$ acts
on $X$ properly and isometrically. Then $\mathcal{A}^{-1}\rtimes G $
is $\KK$-theoretically Poincar\'e dual to
 $\big(\mathcal{A}\hot_{C_0(X)} \ctau (X)\big) \, \rtimes G$, where
 $\mathcal{A}^{-1}$ is the inverse of $\mathcal{A}$ in the Brauer group of
 Morita equivalence classes of continuous trace algebras equipped with a group action.
 We deduce this from a strengthening of Kasparov's
 duality theorem.
As applications we obtain a version of the above Poincar\'e duality
with $X$ replaced by a compact $G$-manifold $M$ and Poincar\'e
dualities for twisted group algebras if the group satisfies some
additional properties related to the Dirac dual-Dirac method for the
Baum-Connes conjecture.
   \end{abstract}


\thanks{This research was supported by the EU-Network \emph{Quantum
  Spaces and Noncommutative Geometry} (Contract HPRN-CT-2002-00280)
  and the \emph{Deutsche Forschungsgemeinschaft} (SFB 478) and by the
  National Science and Engineering Research Council of Canada Discovery
  Grant program.}

\maketitle

\section{Introduction}
\label{sec:intro}

Various versions of Poincar\'e duality in $\KK$-theory and in
equivariant $\KK$-theory played a key role in  the study of the
Baum-Connes and Novikov conjectures and in index theory of
pseudodifferential operators.  In the  past few years, instances of
Poincar\'e duality have appeared, in addition, in quite different
contexts, \emph{e.g.} dynamics (\cite{KP}), and in connection with
Gromov hyperbolic groups (\cite{Em, Em1}). For an extensive
discussion of Poincar\'e duality in the context of operator algebras
we refer to the recent preprint \cite{BMRS} by Brodzki  et al.

The simplest type of duality is the non-equivariant sort, which
was formulated by Connes, Kasparov and Skandalis;
 we refer to this as \emph{noncommutative
duality}. Two C*-algebras $\Lambda$ and $\widehat{\Lambda}$
are noncommutative duals (or
are Poincar\'e dual) if
there is given, for every pair of C*-algebras $(A,B)$, an
isomorphism
$$\Phi_{A,B}\colon \KK(\Lambda\hot A , B) \stackrel{\cong}{\longrightarrow} \KK(A, \widehat{\Lambda}\hot B),$$
and such that $(A,B) \mapsto \Phi_{A,B}$ is natural with respect to
composition in $\KK$ and with respect to external products. Such a
duality is equivalent to the existence of a pair $\widehat{\Delta}
\in \KK(\C, \widehat{\Lambda} \hot \Lambda )$ and $\Delta \in
\KK(\Lambda\hot \widehat{\Lambda} , \C)$ satisfying the `zigzag
equations' of the theory of adjoint functors. One of the simplest
examples is the C*-algebra of a finite group, where duality
is equivalent to a very modest extension of the usual statement of
the Green-Julg theorem, that
$$\KK(C^*(H)\hot A , B) \cong
\KK^H(A, B) \cong \KK(A, C^*(H)\hot B)$$
when $A$ and $B$ are \emph{trivial} $H$-C*-algebras. (The precise
\emph{form} of the isomorphism is of course important.)
The corresponding class $\Delta \in \KK(C^*(H)\hot C^*(H),\C)$ is
that of the $*$-homomorphism
$$a\hot b \mapsto \lambda (a)\rho (b) \in \K\bigl( \ell^2H),$$
where $\lambda$ is the left regular representation and
$\rho$ the right.

The first purpose of this article is to generalize duality for
C*-algebras of finite groups in a very obvious way, namely by
extending to the situation where a discrete, possibly infinite
group $G$, acts properly and co-compactly
on a manifold $X$. In such a situation one can fix a $G$-invariant
Riemannian metric on $X$. We are then in the situation of Kasparov's
duality theorem \cite[Theorem 4.9]{Kas1}, which gives a canonical
isomorphism
\begin{equation}\label{Kas-PD}
\operatorname{RKK}^G(X; A,B) = \mathcal{R}\KK^G\bigl(X; C_0(X,
A), C_0(X, B)\bigr)
  \cong \KK^G\bigl( \ctau (X)\hot A , B)
  \end{equation}
for
arbitrary $G$-$C^{*}$-algebras $A$ and $B$, where
$\Ctau$ is the C*-algebra of $C_0$-sections of the Clifford bundle
 over $X$. It is not very hard to show that Kasparov's duality
implies that the $C^{*}$-algebras
 $C_0(X)\rtimes G$ and $\ctau  (X)\rtimes G$ are Poincar\'e
 dual. It is of interest to have explicit descriptions of the cycles giving rise to the
duality
classes $\Delta$ and $\widehat{\Delta}$, but these will be given elsewhere.

The main contribution of this paper is to extend this result by
stating and proving a twisted version of it, namely that if
$C_0(X,\delta)$ is a continuous trace algebra, i.e., the section
algebra of
  a locally
 trivial $G$-equivariant bundle of elementary C*-algebras representing
 an element $\delta$ of the Brauer group $\mathrm{Br}_G(X)$ (see \cite{CKRW}), then
 $C_0(X,\delta^{-1})\rtimes G $ is Poincar\'e dual to
 $\big(C_0(X,\delta)\hot_{C_0(X)} \ctau (X)\big) \, \rtimes G$, where
 $\delta^{-1}$ is the inverse of $\delta$ in the Brauer group.
To prove this, we need to go back to the argument of the previous
paragraph. What is needed is an extension of
 Kasparov's duality
to the case where $C_0(X, A)$ is replaced by a section algebra $\Gamma_0(E)$
of a smooth $G$-equivariant locally trivial bundle of C*-algebras $p:E\to X$.
 For purposes of this extension, $G$ is allowed to be a
locally compact group and we do not assume the action to be proper,
but we require that the bundle $E$
is endowed with
 a $G$-invariant connection (see Section \ref{sec-twistorb} below for the precise requirements).
 Such connections clearly exist when the bundle $E$ is smooth and a discrete
group $G$ acts smoothly and properly or if a compact group $G$ acts
smoothly on $E$.

We give the proof of our version of Kasparov's duality in Section
\ref{sec-twistorb} below (Theorem \ref{thm-kas-bundle}). In Section
\ref{sec-prop} we first apply this result to obtain Poincar\'e
duality for twisted $K$-theory of compact manifolds equivariant
under an action of a compact group $G$ (Corollary
\ref{cor-compact}). We then obtain our duality result for crossed
products $C_0(X,\delta)\rtimes G$ (Theorem
\ref{thm-poincare-bundle}). As applications, we obtain in
\mbox{Section \ref{sec-twistedgroup}} Poincar\'e dualities for
crossed products $C(M,\delta)\rtimes G$, where $G$ acts
isometrically on the compact manifold $M$ and $\delta$ is a smooth
element of the Brauer group $\Br_G(M)$ and for group algebras
$C^*(G,\om)$ twisted by 2-cocycles $\om\in Z^2(G,\TT)$, if $G$
satisfies some additional properties related to the Dirac-dual Dirac
method for the Baum-Connes conjecture (Theorem \ref{thm-spin} and
Theorem \ref{thm-twist}). A special case of the latter has been
obtained in  \cite[Example 2.6]{BMRS} (see also \cite[Example
2.17]{BMRS}). In Section \ref{sec-nonGcompact}, we extend our
results to a version of Poincar\'e duality on non-$G$-compact
manifolds using a compactly supported $K$-homology theory. Finally,
in the course of proving the main results, we needed to prove the existence of
$G$-equivariant Hermitian connections in certain cases.
However in order to maintain the flow of ideas, we have placed these
auxiliary discussions in the Appendix.


 \begin{remark}\rm
 In the final stages of
 preparation  of this note, a preprint of
 J.-L. Tu has come to our attention (\cite{Tu-P}), in which he proves
 Poincar\'e duality for twisted $K$-theory for compact manifolds
 equivariant under a compact Lie group action, i.e., he gives a
 proof of Corollary \ref{cor-compact} below (with some minor differences
 in the assumptions).
 Although the results have been obtained independently, there
 are noticeable similarities in the method of proof.
\end{remark}

\section{Kasparov's Poincar\'e duality for bundles}\label{sec-twistorb}

If $G$ is a locally compact group, then by a {\em $G$-manifold} we shall always understand
a complete Riemannian manifold $X$ equipped with an action of $G$ by isometric
diffeomorphisms. The section algebra $\Ctau$ of the complex Clifford bundle
of $TX$ (which, by abuse of notation,
we shall also refer to simply as the {\em Clifford bundle of $X$}) is then equipped with a
canonical action of $G$.

Let $p:E\to X$ be a $G$-equivariant locally trivial bundle of
$C^*$-algebras. For such a bundle $p:E\to X$ we denote by
$\Gamma_0(E)$ the algebra of $C_0$-sections of $E$. For simplicity,
we will use the notation $\Gamma_{\tau}(E)$ for
$\Gamma_0(E)\hot_{C_0(X)} C_{\tau}(X)$. We would like to formulate
the following extension of Kasparov's first Poincar\'e duality
(\ref{Kas-PD}) to such bundles:

\begin{equation}\label{PD-bundle}
\RKK^G(X;\Gamma_0(E)\hot A,\,C_0(X,B))\cong
\KK^G(\Gamma_{\tau}(E)\hot A,\,B).
\end{equation}
We refer to \cite[2.19]{Kas1} for the definition of the groups
$\operatorname{RKK}^G(X; A,B)$ for $G$-C*-algebras $A,B$ and
$\RKK^G(X; \mathcal A,\mathcal B)$ for $C_0(X)$-algebras $\mathcal
A, \mathcal B$ and $X$ a locally compact $G$-space. Recall that
$\operatorname{RKK}^G(X; A,B)=\RKK^G(X; C_0(X,A), C_0(X,B))$.

\begin{remark}\label{rem-cex}
We should note that equation (\ref{PD-bundle})
cannot hold for arbitrary $C_0(X)$-algebras, even if $G$ is the trivial group.
To see a simple example  let $X=S^2$  and let $x$ be any fixed element of $S^2$.
Consider $\CC$ as a $C(S^2)$ algebra via the module action
$f\cdot \lambda:=f(x)\lambda$.
Then one easily checks that $\RKK(S^2; \CC, C(S^2))=\{0\}$.
Indeed, if $[\mathcal H, \phi, T]$ is a cycle representing a class in that group,
we may first assume by standard arguments  that  $\phi(\lambda)=\lambda\id_{\mathcal H}$ for
all $\lambda\in \CC$. Then the
Hilbert $C(S^2)$-module $\mathcal H$ consists of continuous sections of a
continuous bundle of Hilbert spaces over $S^2$ such that $f\cdot\xi=f(x)\xi$
for all $f\in C(S^2)$.
But this implies that $f(x)\xi(y)=(f\cdot\xi)(y)=(\xi\cdot f)(y)=\xi(y) f(y)$ for all
$f\in C(S^2)$  and $y\in S^2$. This implies $\xi(y)=0$ for all $y\neq x$, and
then $\xi=0$ by continuity.
On the other hand, we have
$$\KK( C_{\tau}(S^2)\hotimes_{C(S^2)}\CC, \CC)=\KK(\Cl_2,\CC)=\ZZ,$$
where $\Cl_2$ denotes the second complex Clifford algebra.
\end{remark}

Before we proceed, let us investigate in detail the original Kasparov
Poincar\'e duality (\ref{Kas-PD}). We first recall some
definitions and notation of Kasparov from \cite[4.3 and Definition
4.4]{Kas1}.

\begin{notations}\label{r(x)}
Let $X$ be a Riemannian $G$-manifold with  distance function
$\rho:X\times X\to [0,\infty)$. Then  there is a continuous
$G$-invariant positive function $r$ such that for any $x,\,y\in X$
with $\rho(x,y)<r(x)$, there is a unique geodesic segment joining
$x$ and $y$. It can be defined as follows (see \cite[4.3]{Kas1}):
Let $c_0(x)$ be the supremum of absolute values of all sectional
curvatures of $X$ at $x$, and $c_1(x)$ the supremum of $c_0(y)$ for
all $y$ with $\rho(x,y)\leq 1$. Then $r(x)=(c_1(x)+1)^{-1/2}$
satisfies these requirements.

Now let $U_x$ be an open Riemannian ball of radius $r(x)$ around
$x$. Denote by $\Theta_x$ the covector field in $U_x$ which at the
point $y\in U_x$ is given by
$$\Theta_x(y)=\frac{\rho(x,y)}{r(x)}(d_y\rho)(x,y),$$ where $d_y$
means the exterior derivative in $y$. Consider the C*-algebra
$C_{\tau}(U_x)$ as a Hilbert module over $C_{\tau}(X)$. Then
$\Theta_x\in M(C_{\tau}(U_x))$, and ${\Theta_x}^2-1\in C_0(U_x)$.
The field of pairs $\{(C_{\tau}(U_x),\Theta_x)_{x\in X}\}$ defines
an element of the group $\operatorname{RKK}^G(X;\C,C_{\tau}(X))$
which will be denoted by $\Theta_X$. Alternatively, one may consider
$\Theta_X$ as a pair $(\mathcal{F}_U,\Theta)$, where
$$U=\{(x,y)\in X\times X\,|\,\rho(x,y)<r(x)\},\quad
\mathcal{F}_U=\big(C_0(X)\hot C_\tau(X)\big)|_U,$$ and
$\Theta\in\L(\mathcal{F}_U)$ is given by the family $\{\Theta_x\}$
(see \cite[Definition 4.4]{Kas1}). We call $\Theta_X$ the {\em weak
dual-Dirac element of $X$}.
\end{notations}

Kasparov's isomorphism (\ref{Kas-PD}) is given by the map
\begin{equation}\label{mu}
\begin{array}{c}\mu:\operatorname{RKK}^G(X; A,B)\longrightarrow\KK^G\bigl( \ctau (X)\hot A , B)\\
  \mu(\alpha)=\sigma_{X,C_{\tau}(X)}(\alpha)\hot_{C_{\tau}(X)}d_X
\end{array}
\end{equation}
with its inverse given by the map
\begin{equation}\label{nu}
\begin{array}{c}\nu:\KK^G\bigl( \ctau (X)\hot A , B)\longrightarrow\operatorname{RKK}^G(X; A,B)\\
  \nu(\beta)=\Theta_X\hot_{C_0(X)\hot
  C_{\tau}(X)}\sigma_{C_0(X)}(\beta),
\end{array}
\end{equation}
(see \cite[Theorem 4.9]{Kas1}) where $d_X$ is the Dirac class defined in \cite[4.2]{Kas1}.
Now, in a first attempt to prove
the isomorphism (\ref{PD-bundle}) one can see that the formula from
(\ref{mu}) still gives us a well-defined homomorphism
\begin{equation}\label{mu1}
\begin{array}{c}\mu:\RKK^G(X;\Gamma_0(E)\hot
A,\,C_0(X,B))\longrightarrow\KK^G(\Gamma_{\tau}(E)\hot A,\,B)\\
\mu(\alpha)=\sigma_{X,C_{\tau}(X)}(\alpha)\hot_{C_{\tau}(X)}d_X,
\end{array}
\end{equation}
even if we would replace $\Gamma_0(E)$ by any $C_0(X)$-algebra.
But, unfortunately, the inverse map $\nu$ of (\ref{nu}) has no obvious extension to the
bundle case, and it is the main technical result of this paper to show that  in
many important situations one can find a substitute for (\ref{nu}) which will do the job
(it follows from Remark \ref{rem-cex} above that a similar result would be impossible for
general $C_0(X)$-algebras).
As we shall see below, the main ingredient in this construction will be the construction
of a  certain element  $\Theta_{E}\in \RKK^G(X;\Gamma_0(E),C_0(X)\hot\Gamma_{\tau}(E))$
(with $C_0(X)$-action on the first factor of the algebra
$C_0(X)\hot\Gamma_{\tau}(E)$), which will replace the element $\Theta_X$
in (\ref{nu}).

\begin{remark}\label{rem-r(x)}
{\bf (1)} Notice that the class $\Theta_X\in
\operatorname{RKK}^G(X;\C,C_{\tau}(X))$ does not depend on the
particular choice of the function $r:X\to (0,\infty)$. In fact, if
$r_0, r_1:X\to (0,\infty)$ are two such maps, then the elements
$\Theta_t\in \operatorname{RKK}^G(X;\C,C_{\tau}(X))$ constructed as
above via the distance functions $r_t(x):=tr_1(x)+(1-t)r_0(x)$ give
a homotopy connecting the elements corresponding to $r_0,r_1$. In
particular, if $x\mapsto r(x)$ has a lower bound $\eta>0$ (e.g., if
$G\backslash X$ is compact), then we may replace the function $r$ by
the constant $\eta$ in the construction of $\Theta_X$. In that case,
the set $U$ is clearly symmetric.
\\
{\bf (2)} The pair $(\mathcal{F}_U,\Theta)$ also determines a class
$$\Theta_U\in \RKK^G(X; C_0(X), (C_0(X)\hot C_\tau(X))|_U),$$
 so that the class
$\Theta_X\in\operatorname{RKK}^G(X; \CC, C_{\tau}(X))=\RKK^G(X;
C_0(X), C_0(X)\hot C_\tau(X))$ satisfies the equation
\begin{equation}\label{eq-ThetaU}
\Theta_X=\iota_*(\Theta_U),
\end{equation}
where $\iota: (C_0(X)\hot C_\tau(X))|_U\to C_0(X)\hot C_\tau(X)$
denotes the inclusion map. In what follows below, we shall often
work with the element $\Theta_U$.
\end{remark}

\begin{remark}\label{rem-ThetaE-idea} The first idea for the construction of the
element $\Theta_E$ was to
consider the class $\sigma_{X,\Gamma_0(E)}(\Theta_X)\in
\RKK^G(X;\Gamma_0(E),\Gamma_0(E)\hot C_{\tau}(X))$, but
it turns out that the algebra in the second variable of $\RKK^G(X;\cdot,\cdot)$
is not quite right.
We probably could define $\Theta_E$
via $\sigma_{X,\Gamma_0(E)}(\Theta_X)$ if we had a $G$-equivariant
$C_0(X\times X)$-linear isomorphism
\begin{equation}\label{C(XX)-linear} \Gamma_0(E)\hot C_{\tau}(X)\cong
C_0(X)\hot\Gamma_{\tau}(E).
\end{equation}
Note that the two algebras are obtained by pulling back the bundle
$p:E\rightarrow X$ via the projection maps $p_0,p_1:X\times
X\rightarrow X$ on the first and second variable,
respectively, and then by tensoring $C_{\tau}(X)$ on the second
variable of the section algebras of the pull-backs. We shall define
$\Theta_E$ in Definition \ref{defn-ThetaE} based on this idea after
we discus some further preliminaries.
\end{remark}

\begin{notations}\label{P} Recall the construction of the open subset $U$ of $X\times X$ in Notation
\ref{r(x)}. The important property of $U$ is that for each $(x,y)\in
U$ there exists a unique geodesic $\gamma_{x,y}:[0,1]\to X$ joining
$x$ and $y$. Then the projection on the first variable $p_0:U\to X$
and the projection on the second variable $p_1:U\to X$ are
$G$-equivariantly homotopic given by the homotopy
 $$p_t:U\longrightarrow X;\,p_t(x,y)=\gamma_{x,y}(t).$$
Let $P:U\times[0,1] \to X, \;P(x,y,t)=p_t(x,y)\quad\text{and}\quad
P_t: U\times [0,1] \to X;\; P_t(x,y,s)=p_t(x,y).$ Then, forgetting the
$G$-action, it is a well-known fact from the theory of fibre bundles
(e.g. see \cite{Huse}) that there exists
 a bundle isomorphism
\begin{equation}\label{bundle-iso}
\varphi: P^*E\to P_0^*E\cong [0,1]\times p_0^*E \end{equation}
which fixes the base $U\times [0,1]$ of these bundles, where, as usual,
$P^*E$ and $P_0^*E$ denote the pull-backs of the bundle $E$ via the
maps $P$ and $P_0$.
\end{notations}

\begin{remark}\label{rem-feasible}
It is not clear whether or not one can always arrange for
the map $\varphi$ in (\ref{bundle-iso}) to be
$G$-equivariant.  If it {\em is} possible,  we shall call $p:E\to X$
a {\em feasible} $G$-bundle over $X$.

We certainly have feasibility if $E=X\times A$ is a trivial bundle,
$A$ is a $G$-C*-algebra and the action on $E$ is given by the
diagonal action, in which case the isomorphism (\ref{PD-bundle}) is
equivalent to Kasparov's original result (\ref{Kas-PD}). Moreover,
it is shown in \cite{AS} and \cite[Section 3.2]{Tu-P} that $p:E\to
X$ is always $G$-isomorphic to a feasible bundle in our sense if the
fibres of $E$ are the compact operators $\K$, $G$ is a compact Lie
group and the bundle $p:E\to X$ satisfies a certain local condition
which is spelled  out in \cite[\S 6]{AS}. If $G$ is a compact group
which acts {\em trivially} on the base $X$ of a locally trivial
bundle of compact operators $p:E\to X$  one can use the
classification of such bundles as given in \cite{EW1} to show that
they are always feasible at least if we allow to pass to a Morita
equivalent bundle (which will not affect equation
(\ref{PD-bundle})). However, the details of this could easily fill
one or two pages, and we do not think that this case is interesting
enough to justify  this effort.
\end{remark}

In general, obtaining  feasibility can be a non-trivial problem.
In the Appendix (see Lemma \ref{lem-connection} and Proposition \ref{lem-trans} below)
we shall prove:

\begin{prop}\label{prop-feasible}
Suppose that $p:E\to X$ is a locally trivial  smooth
 $G$-bundle of C*-algebras over the $G$-manifold $X$ (this implies in particular, that
 the action of $G$ is by smooth bundle automorphisms).
 Then $p:E\to X$
is feasible if there exists a $G$-invariant Hermitian connection on $E$. In particular,
if $G$ is compact or if $G$ is discrete and acts properly on $X$, then $p:E\to X$
is feasible.
\end{prop}

\begin{remark}\label{rem-cliff-con} Notice that for every $G$-manifold
$X$ the Levi-Civita connection induces a canonical $G$-invariant
connection on the Clifford bundle, so the Clifford bundle of a $G$-manifold $X$
is always feasible.
\end{remark}

\begin{remark}\label{wockel} A recent paper by M\"uller and Wockel (\cite{MW}) actually
implies that every locally trivial bundle $p:E\to X$ of compact operators over a manifold
$X$ is isomorphic to a smooth bundle. So the bundles of our main interest, namely
those whose section algebras are  continuous trace algebras can always
be chosen as to be smooth.
 But, unfortunately, it is not clear whether any
(proper) group action on the given algebra can be realized, up to
Morita equivalence, by a smooth action on this bundle. So it is not
clear whether our smoothness assumptions can be avoided in that
situation. Note that it is definitely not true that every
automorphism of $A$ comes from a smooth automorphism of $E$, even if
the bundle is trivial and the group acts trivially on the base. The
reason is that strongly continuous maps from $X$ to $\PU$ are not
automatically smooth.
\end{remark}

Suppose for now that $\varphi:P^*E\to P_0^*E$ is a $G$-equivariant bundle
isomorphism. Then it induces a $G$-equivariant and
$C_0(U\times[0,1])$-linear isomorphism of the section algebras
\begin{equation}\label{phi}
\phi: \Gamma_0(P^*E)\stackrel{\cong}{\to} \Gamma_0(P_0^*E)
\end{equation}
which restricts on the sets $U\times\{t\}\cong U$ to  $G$-equivariant and $C_0(U)$-linear
isomorphisms
\begin{equation}\label{phi_t}
\phi_t: \Gamma_0(p_t^*E)\stackrel{\cong}{\to}  \Gamma_0(p_0^*E)
\end{equation}
for all $t\in[0,1]$ (since $\varphi$ restricts  to isomorphisms $\varphi_t:p_t^*E\to p_0^*E$
on $U\times\{t\}\cong U$).
Note further that
\begin{equation}\label{eq-p_0E}
\Gamma_0(p_0^*E)=(\Gamma_0(E)\hot C_0(X))|_U
\quad\text{and}\quad
\Gamma_0(p_1^*E)=(C_0(X)\hot
\Gamma_0(E))|_U,
\end{equation}
where for any continuous C*-algebra bundle $E$ over some base space
$Y$ and any locally compact subset $V$ of $Y$  we put
$\Gamma_0(E)|_V:=\Gamma_0(E|_V)$ (note that $\Gamma_0(E)|_V$ can also be realized
as the balanced tensor product
$C_0(V)\hot_{C_0(Y)}\Gamma_0(E)$ (e.g. see \cite{RW0}) which carries
canonical $C_0(V)$- {\bf and}  $C_0(Y)$-algebra structures). Recall
from (\ref{phi_t}) that we have a $C_0(U)$-linear (and hence also
$C_0(X\times X)$-linear) isomorphism
\begin{equation}\label{eq-phi1}
\phi_1: (\Gamma_0(E)\hot C_0(X))|_U\to (C_0(X)\hot \Gamma_0(E))|_U.
\end{equation}
If we consider both algebras as $C_0(X)$-algebras with respect to
the second component,  it follows from a simple argument of
associativity of balanced tensor product (see Lemma \ref{tensor
product} below) that we obtain a $G$-equivariant $C_0(X)\hot
C_\tau(X)$-linear isomorphism
\begin{equation}\label{eq-barphi}
\overline{\phi}_1:=\phi_1\hot_{C_0(X)}\id_{C_{\tau}(X)}: (\Gamma_0(E)\hot C_\tau(X))|_U
\stackrel{\cong}{\to} (C_0(X)\hot \Gamma_\tau(E))|_U.
\end{equation}
\begin{lem}\label{tensor product}
Assume that $A$ is a $C_0(X\times Y)$-algebra, $B$ is a
$C_0(X)$-algebra and $C$ is a $C_0(Y)$-algebra. Then
$B\hot_{C_0(X)}A$ and $A\hot_{C_0(Y)}C$ are $C_0(X\times
Y)$-algebras and there is a canonical isomorphism
$$\left(B\hot_{C_0(X)}A\right)\hot_{C_0(Y)}C\cong
B\hot_{C_0(X)}\left(A\hot_{C_0(Y)}C\right)$$ as $C_0(X\times
Y)$-algebras. This isomorphism is also $C_0(X\times
Y)\hot_{C_0(Y)}C\cong C_0(X)\hot C$-linear.
\end{lem}
\begin{proof}
It is clear that the $C_0(X\times Y)$-structure on $A$ induces
canonical $C_0(X\times Y)$-structures on $B\hot_{C_0(X)}A$ and
$A\hot_{C_0(Y)}C$. It is then easy to check that the map
$$(b\hot a)\hot c\mapsto
b\hot(a\hot c)$$ satisfies all requirements of the lemma.
\end{proof}

\begin{defn}\label{defn-ThetaE} For a feasible bundle
$p:E\longrightarrow X$, we define $\Theta_E$ as the element
$$\Theta_E= \iota_*\big(\overline{\phi}^{-1}_{1,*}\big(\sigma_{X,
\Gamma_0(E)}(\Theta_U)\big)\big) \in\RKK^G(X; \Gamma_0(E),
C_0(X)\hot \Gamma_\tau(E)),$$ where $\Theta_U$ is from Remark
\ref{rem-r(x)}(2). Here $\sigma_{X,\Gamma_0(E)}(\Theta_U)\in
\RKK^G(X; \Gamma_0(E), (\Gamma_0(E)\hot  C_{\tau}(X))|_U)$ is the
element obtained by tensoring $\Theta_U$ over $C_0(X)$ with
$\Gamma_0(E)$,
$$\overline{\phi}_1^{-1}:(\Gamma_0(E)\hot  C_{\tau}(X))|_U\to (C_0(X)\hot \Gamma_\tau(E))|_U$$
is the inverse of the isomorphism of (\ref{eq-barphi}), and
$\iota: (C_0(X)\hot \Gamma_\tau(E))|_U\to C_0(X)\hot \Gamma_{\tau}(E)$ denotes
inclusion.
\end{defn}

\begin{remark}\label{rem-ThetaE} For later use it is necessary to give a precise description of
a Kasparov triple  corresponding to $\Theta_E$. For this it is
convenient to introduce some further notation. For $t\in [0,1]$
recall that $p_t:U\to X$ is defined by $p_t(x,y)=\gamma_{x,y}(t)$,
with $\gamma_{x,y}$ the geodesic joining $x$ with $y$. We then have
an obvious equation
\begin{equation}\label{eq-FU}
\F_U\hot_{C_0(X),p_t}\Gamma_0(E):= \F_U\hot_{C_0(U)} \Gamma_0(p_t^*E)
\end{equation}
where $\F_U\hot_{C_0(X),p_t}\Gamma_0(E)$ is the balanced tensor product of
$\F_U=\big(C_0(X)\hot C_\tau(X)\big)|_U$ with $\Gamma_0(E)$ when $\F_U$ carries the
$C_0(X)$ structure induced by the  map $p_t:U\to X$.
Note that, in particular, we have
$$\F_U\hot_{C_0(X),p_0}\Gamma_0(E)\cong \big(\Gamma_0(E)\hot C_\tau(X)\big)|_U
\quad\text{and}\quad
\F_U\hot_{C_0(X),p_1}\Gamma_0(E)\cong \big(C_0(X)\hot \Gamma_\tau(E)\big)|_U.$$
Following Definition \ref{defn-ThetaE}, we obtain
\begin{equation}\label{eq-ThetaE}
\Theta_E=[(\F_U\hot_{C_0(X),p_1}\Gamma_0(E), \psi_E, \Theta\hot 1)],
\end{equation}
with left action
of $\Gamma_0(E)$ on this module given by the formula
\begin{equation}\label{eq-psiE}
\psi_{E}(f)\cdot \xi := (\overline{\phi}_1)^{-1}\left(f\cdot
\overline{\phi}_1(\xi)\right)
\end{equation}
where the product $f\cdot \overline{\phi}_1(\xi)$ is given via the canonical action of $\Gamma_0(E)$
on
$$\F_U\hot_{C_0(X),p_0}\Gamma_0(E)
\stackrel{\overline{\phi}_1^{-1}}{\cong}
\F_U\hot_{C_0(X),p_1}\Gamma_0(E).$$ Since this isomorphism is
$\F_U$-linear (\ref{eq-phi1}),
 it follows that the operator in $\Theta_E$ is given by Kasparov's operator $\Theta$
acting on the first variable of $\F_U\hot_{C_0(X),p_1}\Gamma_0(E)$.
We now state the main result of this section:

\begin{thm}\label{thm-kas-bundle} Let
$p:E\longrightarrow X$ be a feasible locally trivial $G$-equivariant
$C^*$-algebra bundle over a $G$-manifold $X$ and recall that
$\Gamma_{\tau}(E)= \Gamma_0(E)\hot_{C_0(X)}C_{\tau}(X)$. Then, for any pair of
$G$-algebras $A$ and $B$, we have a natural isomorphism
\begin{align*}
\mu:\RKK^G(X;\Gamma_0(E)&\hot A,\,C_0(X,B))\longrightarrow
\KK^G(\Gamma_{\tau}(E)\hot A,\,B)\\
&\mu(\alpha):= \sigma_{X,C_{\tau}(X)}(\alpha)\hot_{C_{\tau}(X)}d_X
\end{align*}
with inverse given by
\begin{align*}\nu:\KK^G(\Gamma_{\tau}(E)&\hot
A,\,B)\longrightarrow \RKK^G(X;\Gamma_0(E)\hot
A,\,C_0(X,B))\\
&\nu(\beta):=\Theta_E\hot_{C_0(X)\hot
\Gamma_{\tau}(E)}\sigma_{C_0(X)}(\beta),
\end{align*}
where $\Theta_E$ is the class as defined in Definition
\ref{defn-ThetaE}.
\end{thm}

In Section 4 we shall prove an analogue of the above theorem
for a manifold with boundary (see Theorem \ref{thm-kas-bundle1} below).
The main work for the proof of the Theorem \ref{thm-kas-bundle}
will be done in the proofs of the technical lemmas, Lemma  \ref{lem4.5}
and Lemma \ref{lem4.6} below.
For later use in the proof of Lemma \ref{lem4.5}, we
actually have to extend the construction of $\Theta_E$ as follows:
Consider $P:U\times[0,1]\to X; P(x,y,t)=p_t(x,y)$. We want to
construct a certain  element
$$\Theta_P\in \RKK^G\big(X\times[0,1]; \Gamma_0(E)[0,1], \F_U[0,1]\hot_{C_0(U\times[0,1])}\Gamma_0(P^*E)\big),$$
where for any algebra $A$ we write $A[0,1]$ for $C\big([0,1],
A\big)$.  Since $E$ is feasible we have
an isomorphism
$$\phi: \Gamma_0(P^*E)\to \Gamma_0(P_0^*(E))\cong \Gamma_0(p_0^*E)[0,1]$$
which then induces an $\F_U[0,1]$-linear isomorphism
\begin{equation}\label{eq-barphi1}
\overline{\phi}: \F_U[0,1]\hot_{C_0(U\times[0,1]) }\Gamma_0(P^*E)\stackrel{\cong}{\to}
\F_U[0,1]\hot_{C_0(U\times[0,1]) }\Gamma_0(P_0^*E).
\end{equation}
Moreover, we have an obvious isomorphisms
\begin{equation}\label{eq-isoP}
\F_U[0,1]\hot_{C_0(U\times[0,1]) }\Gamma_0(P_0^*E)
\cong \big(\F_U\hot_{C_0(X),p_0}\Gamma_0(E)\big)[0,1].
\end{equation}
Thus we can  define the element $\Theta_P$ by
\begin{equation}\label{eq-ThetaP}
\Theta_P:=\overline{\phi}^{-1}_*\left(\sigma_{C[0,1]}\big(\sigma_{X,\Gamma_0(E)}(\Theta_U)\big)\right).
\end{equation}
Evaluation of this element at $t\in [0,1]$ gives an element
$$\Theta_t\in \RKK^G(X; \Gamma_0(E), \F_U\hot_{C_0(X), p_t}\Gamma_0(E))$$
which is given by the Kasparov triple
\begin{equation}\label{eq-phit}
\Theta_t=[(\F_U\hot_{C_0(X), p_t}\Gamma_0(E), \psi_t, \Theta\hot 1)]
\end{equation}
with $\psi_t:\Gamma_0(E)\to \L(\F_U\hot_{C_0(X), p_t}\Gamma_0(E))$ given by
the formula
\begin{equation}\label{eq-psit}
\psi_t(f)\xi:=(\overline{\phi}_t)^{-1}\left(f\cdot
\overline{\phi}_t(\xi)\right)
\end{equation}
where $\overline{\phi}_t: \F_U\hot_{C_0(U)}\Gamma_0(p_t^*E)\to
 \F_U\hot_{C_0(U)}\Gamma_0(p_0^*E)$
 is the isomorphism induced by the isomorphism $\phi_t:\Gamma_0(p_t^*E)\to \Gamma_0(p_0^*E)$
of (\ref{eq-phit}).
In particular, it  follows from our constructions that
\begin{equation}\label{eq-Thetaglobal}
\iota_*(\Theta_1)=\Theta_E\quad\text{and}\quad
\Theta_0=\sigma_{X,\Gamma_0(E)}(\Theta_U).
\end{equation}
 \end{remark}

The following lemma is an extension of \cite[Lemma 4.5]{Kas1}:

\begin{lem}\label{lem4.5}
Let $A$ and $B$ be $G$-algebras.
For any $\alpha\in \RKK^G(X;\Gamma_0(E)\hot
A,\,C_0(X,B))$, the equation
\begin{equation}\label{alpha}
\alpha\hot_{C_0(X)}\Theta_X=
\Theta_E\hot_{C_0(X)\hot\Gamma_{\tau}(E)}\left(\sigma_{C_0(X)}\big(\sigma_{X,C_{\tau}(X)}\left(\alpha\right)\big)\right).
\end{equation}
holds in $\RKK^G(X;\Gamma_0(E)\hot
A,\,C_0(X)\hot C_{\tau}(X)\hot B)$.

\end{lem}

\begin{proof}
Since taking Kasparov products over $C_0(X)$ in $\RKK^G(X;\cdot,\cdot)$-theory is commutative
by \cite[Proposition 2.21]{Kas1}, we get
$$\alpha\hot_{C_0(X)}\Theta_X=\Theta_X\hot_{C_0(X)}\alpha=\sigma_{X,\Gamma_0(E)}(\Theta_X)\hot_{\Gamma_0(E)}\alpha.$$
Hence it is enough to show that
$$\sigma_{X,\Gamma_0(E)}(\Theta_X)\hot_{\Gamma_0(E)}\alpha=
\Theta_E\hot_{C_0(X)\hot\Gamma_{\tau}(E)}\left(\sigma_{C_0(X)}\big(\sigma_{X,C_{\tau}(X)}\left(\alpha)\big)\right)\right).$$
Consider the map $P:U\times[0,1]\rightarrow X, (x,y,t)\mapsto p_t(x,y)$
and  the pull-back
$$P^*(\alpha)\in \RKK^G\left(U\times[0,1];\,\Gamma_0(P^*E)\hot A,\,C_0(U\times[0,1])\hot B\right)$$
as in \cite[Proposition 2.20]{Kas1}.
Let
$$
\begin{array}{rcl}
\lefteqn{\beta:=j_*\left(\sigma_{U\times[0,1],\mathcal{F}_U[0,1]}(P^*\alpha)\hot[\iota]\right)}\\
&\in&\RKK^G\left(X\times[0,1];\,\mathcal{F}_U[0,1]\hot_{C_0(U\times[0,1])}\Gamma_0(P^*E)\hot
A,\,C_0(X)\hot C_{\tau}(X)\hot B[0,1]\right),
\end{array}
$$
where $\iota:\mathcal{F}_U[0,1]\hot B\hookrightarrow
C([0,1])\hot C_0(X)\hot C_{\tau}(X)\hot B$ denotes inclusion and $j_*$ is the
map which changes the  $C_0(U\times[0,1])$-structure to the
$C_0(X\times[0,1])$-structure by letting $C_0(X)$ act on the first variable in $C_0(U)$. Then the
restrictions of $\beta$ at $0$ and $1$ are
\begin{equation}\label{beta01}
\begin{array}{rcl}
\beta_0&=&[\iota^0]\hot_{C_{\tau}(X)\hot\Gamma_0(E)} \big(\sigma_{C_{\tau}(X)}(\alpha)\big)\\
\beta_1&=&[\iota^1]\hot_{C_0(X)\hot \Gamma_{\tau}(E)}
\left(\sigma_{C_0(X)}\big(\sigma_{X,C_{\tau}(X)}(\alpha)\big)\right),
\end{array}
\end{equation}
where $\iota^0,\iota^1$ denote the inclusions
$$\iota^0:\FUpzero \hot A\hookrightarrow
C_{\tau}(X)\hot\Gamma_0(E)\hot A$$ and
$$\iota^1:\FUpone
\hot A\hookrightarrow C_0(X)\hot \Gamma_{\tau}(E)\hot A.$$
Indeed, if $\alpha$ is represented by the triple
$(\E, \varphi_{\alpha}, T)$,
the isomorphism
$$(C_0(X)\hot
C_{\tau}(X))\hot_{C_0(X),p_0}\E\cong \E\hot C_{\tau}(X)$$
restricted to $U$ gives the isomorphism
$\mathcal{F}_U\hot_{C_0(X),p_0}\E\cong (\E\hot C_\tau(X))|_U$ and then
$$\begin{array}{rcl}
\beta_0&=&\sigma_{U,\mathcal{F}_U}\left(p_0^*(\alpha)\right)\hot[\iota]\\
&=&\left[\left(\mathcal{F}_U\hot_{C_0(U)}(C_0(U)\hot_{C_0(X),p_0}\E),\,\id_{\mathcal{F}_U}\hot(\id_{C_0(U)}\hot
\varphi_{\alpha}),\,1\hot T\right)\right]\\
&=&\left[ \left((\E\hot C_{\tau}(X))|_U,\,\varphi_{\alpha}\hot \id_{C_{\tau}(X)},\,T\hot
1\right)\right].
\end{array}
$$
which is easily seen to coincide with
$[\iota^0]\hot_{C_{\tau}(X)\hot\Gamma_0(E)} \big(\sigma_{C_{\tau}(X)}(\alpha)\big)$.
Similarly,
\begin{align*}
\beta_1&=\sigma_{U,\mathcal{F}_U}(p_1^*(\alpha))\hot[\iota]\\
&=\left[\left(\mathcal{F}_U\hot_{C_0(X),p_1}\E,\,\id_{\mathcal{F}_U}\hot
\varphi_{\alpha},\,1\hot T\right)\right]
\end{align*}
and
\begin{align*}
\lefteqn{[\iota^1]\hot_{C_0(X)\hot \Gamma_{\tau}(E)\hot A}
\left(\sigma_{C_0(X)}\left(\sigma_{X,C_{\tau}(X)}(\alpha)\right)\right)}\\
&=[\iota^1]\hot_{C_0(X)\hot
\Gamma_{\tau}(E)\hot A}
\left[ \left(C_0(X)\hot
(C_{\tau}(X)\hot_{C_0(X)}\E),\,\id_{C_0(X)}\hot
(\id_{C_{\tau}(X)}\hot\varphi_{\alpha}),\,1\hot T\right)\right]\\
&=\left[\left(\mathcal{F}_U\hot_{C_0(X),p_1}\E,\,\id_{\mathcal{F}_U}\hot
\varphi_{\alpha},\,1\hot T\right)\right].
\end{align*}
To complete the proof of the lemma, we consider the element
$$
\Theta_\alpha
:=\sigma_A(\Theta_P)\hot_{\mathcal{F}_U[0,1]\hot_{C_0(U\times[0,1])}(\Gamma_0(P^*E)\hot
A)}\,\,\beta
$$
in $\RKK^G\left(X\times[0,1];\,\Gamma_0(E)\hot
A[0,1],\,C_0(X)\hot C_{\tau}(X)\hot B[0,1]\right)$.
Then by (\ref{eq-Thetaglobal}) and (\ref{beta01}), the evaluations
of $\Theta_{\alpha}$ at $0$ and $1$ are the following:

\begin{align*}
(\Theta_{\alpha})_0&
=\sigma_A\left(\Theta_0\right)\hot_{\mathcal{F}_U\hot_{C_0(U)}(\Gamma_0(p_0^*E)\hot
A)}\,\,\beta_0\\
&=\left(\sigma_A\big(\sigma_{X,\Gamma_0(E)}({\Theta_U})\big)\hot_{\mathcal{F}_U\hot_{C_0(U)}(\Gamma_0(p_0^*E)\hot
A)}\,\,[\iota^0]\right)\hot_{C_{\tau}(X)\hot\Gamma_0(E)\hot
A}\big(\sigma_{C_{\tau}(X)}(\alpha)\big)\\
&=\sigma_{X,\Gamma_0(E)}(\Theta_X)\hot_{\Gamma_0(E)}\alpha
\end{align*}
and
\begin{align*}
(\Theta_{\alpha})_1&
=\sigma_A\left(\Theta_1\right)\hot_{\mathcal{F}_U\hot_{C_0(U)}(\Gamma_0(p_1^*E)\hot
A)}\,\,\beta_0\\
&=\left(\sigma_A({\Theta_1})\hot_{\FUpone \hot
A}\,\,[\iota^1]\right)\hot_{C_{0}(X)\hot\Gamma_{\tau}(E)\hot
A}\,\left(\sigma_{C_0(X)}\big(\sigma_{X,C_{\tau}(X)}(\alpha)\big)\right)\\
&=\sigma_A(\iota^1_*(\Theta_1))\hot_{C_0(X)\hot \Gamma_{\tau}(E)\hot
A}\,\left(\sigma_{C_0(X)}\big(\sigma_{X,C_{\tau}(X)}(\alpha)\big)\right)\\
&=\sigma_A(\Theta_E)\hot_{C_0(X)\hot \Gamma_{\tau}(E)\hot
A}\,\left(\sigma_{C_0(X)}\big(\sigma_{X,C_{\tau}(X)}(\alpha)\big)\right).
\end{align*}
This proves that
$\sigma_{X,\Gamma_0(E)}(\Theta_X)\hot_{\Gamma_0(E)}\alpha$ and
$\sigma_A(\Theta_E)\hot_{C_0(X)\hot \Gamma_{\tau}(E)\hot
A}\,\left(\sigma_{C_0(X)}\big(\sigma_{X,C_{\tau}(X)}(\alpha)\big)\right)$ are homotopic.
\end{proof}

Another important step for the proof of Theorem \ref{thm-kas-bundle} is the following
bundle-version of \cite[Lemma 4.6]{Kas1}.

\begin{lem}\label{lem4.6}
Let
$\Sigma:C_{\tau}(X)\hot\Gamma_{\tau}(E)\longrightarrow\Gamma_{\tau}(E)\hot
C_{\tau}(X)$ denote the canonical flip-isomorphism. Then
\begin{equation}\label{equality}\sigma_{X,\Gamma_{\tau}(E)}\left(\Theta_X\right)=\Sigma_*\left(\sigma_{X,C_{\tau}(X)}\left(\Theta_E\right)\right)
\end{equation}
in $\KK^G\left(\Gamma_{\tau}(E),\Gamma_{\tau}(E)\hot
C_{\tau}(X)\right).$
\end{lem}
\begin{proof} Recall the construction of $\Theta_E$ from Definition
\ref{defn-ThetaE} and Remark \ref{rem-ThetaE}:
$$\Theta_E=[(\FUpone,\,\psi_{E},\Theta)].$$ The fiberwise
action by $\psi_E$ is as follows:
\begin{equation*}
\left(\psi_E(\eta)(g\hot
\xi)\right)(x,y)=g(x,y)\hot\left(\varphi_{x,y,1}\right)^{-1}\left(\eta(x)\right)\xi(y),
\end{equation*}
where
$\eta\in\Gamma_0(E),\,g\hot\xi\in\mathcal{F}_U\hot_{C_0(X),p_1}\Gamma_0(E)$
and
$$\varphi_{x,y,1}: (P^*E)_{(x,y,1)}=E_y\to E_x=(P_0^*E)_{(x,y,1)}$$
is the fibre map of the bundle isomorphism $\varphi: P^*E\to P_0^*E$ of Notations \ref{P}.

>From the above, we have
$$\sigma_{X,C_{\tau}(X)}(\Theta_E)=\left(C_{\tau}(U)\hot_{C_0(X),p_1}\Gamma_0(E),\,
\widetilde{\psi}_1,\,\Theta\hot
1_{\Gamma_0(E)}\right),$$ where
$\widetilde{\psi}_1=\id_{C_{\tau}(X)}\hot \psi_E$ and
$\Theta(x,y)=\frac{\rho(x,y)}{r(x)}\left(d_y\rho\right)(x,y)$ from
Notations \ref{r(x)}. Then the right hand side of (\ref{equality})
is written as follows:
\begin{equation}\label{homotopy-1}
\Sigma_*\left(\sigma_{X,C_{\tau}(X)}(\Theta_E)\right)=\left(C_{\tau}(\widetilde{U})\hot_{C_0(X),p_0}\Gamma_0(E),\,{\psi}^{\Sigma},\,\Theta^{\Sigma}\hot
1\right),
\end{equation}
where $\widetilde{U}=\{(x,y): (y,x)\in U\}$ and
\begin{equation}\label{eq-psiSigma}
\left(\psi^{\Sigma}(f\hot_{C_0(X)}\eta)(g\hot\xi)\right)(x,y)=f(y)\cdot g(x,y)\hot(\varphi_{y,x,1})^{-1}\left(\eta(y)\right)\xi(x)
\end{equation}
for $f\hot_{C_0(X)}\eta\in
C_{\tau}(X)\hot_{C_0(X)}\Gamma_0(E),\,g\hot\xi\in
C_{\tau}(\widetilde{U})\hot_{C_0(X),p_0}\Gamma_0(E)$, and
\begin{equation}\label{eq-ThetaSigma}
\Theta^{\Sigma}(x,y)=\frac{\rho(x,y)}{r(y)}\left(d_x\rho\right)(x,y).
\end{equation}
On the other hand, the element $\sigma_{X,\Gamma_{\tau}(E)}$ is easily seen to be represented by
the triple
\begin{equation}\label{eq-left}
\left(C_{\tau}({U})\hot_{C_0(X),p_0}\Gamma_0(E),\,{\psi}_0,\,\Theta\hot 1\right),
\end{equation}
with  $\Theta$ as above and with
\begin{equation}\label{eq-psi0}
\psi_0\left((f\hot_{C_0(X)}\eta)(g\hot\xi)\right)(x,y)=f(x)\cdot g(x,y)\hot \eta(x)\xi(x),
\end{equation}
for $f\hot_{C_0(X)}\eta\in
C_{\tau}(X)\hot_{C_0(X)}\Gamma_0(E),\,g\hot\xi\in
C_{\tau}({U})\hot_{C_0(X),p_0}\Gamma_0(E)$.

We are going to use Kasparov's homotopy
$\{\left[(C_\tau(U), \mu_t, \Theta_t)\right]\}_{t\in [0,1]}$ between
$\sigma_{X, C_{\tau}(X)}(\Theta_X)$ and the class
$\widetilde{\Theta}\in \KK^G(C_\tau(X), C_\tau(X)\hot C_\tau(X))$
given by the triple
$(C_{\tau}(U), \mu_1,\Theta_1)$, where $C_\tau(U)=(C_\tau(X)\hot C_\tau(X))|_U$,
$\mu_1: C_\tau(X)\to \L(C_\tau(U))$ is given by multiplication on the second factor
and
\begin{equation}\label{eq-theta1}
\Theta_1(x,y)=\frac{\rho(x,y)}{r(x)}\left(d_x\rho\right)(x,y),
\end{equation}
which is carefully described in the proof of \cite[Lemma 4.6]{Kas1}.
We may consider the triple $(C_\tau(U), \mu_t, \Theta_t)$ as  a cycle for a class in
$\KK^G(C_\tau(X),C_{\tau}(U))$ (note that the homomorphism $\mu_t: C_{\tau}(X)\to \L(C_{\tau}(U))$
is denoted $\varphi_t$ in \cite{Kas1}). One can easily check that it determines a class
in $\RKK^G(X, C_{\tau}(X), C_{\tau}(U))$, when the $C_0(X)$-structure on $C_\tau(U)$ is given by
the formula
$$f\cdot \xi (x,y):=f(p_t(x,y))\xi(x,y)\quad\text{for all}\; f\in C_0(X), \xi\in C_{\tau}(U).$$
Thus, we see that the family $\{\left[(C_\tau(U), \mu_t, \Theta_t)\right]\}_{t\in [0,1]}$ determines a class
$$\widetilde{\Theta}\in \RKK^G(X, C_\tau(X), C_{\tau}(U)[0,1])$$
with respect to the $C_0(X)$-structure on $C_\tau(U)[0,1]$
given by
$$f\cdot \eta(x,y,t)=f(p_t(x,y))\eta(x,y,t)\quad\text{for all}\; f\in C_0(X), \eta\in C_{\tau}(U)[0,1].$$
Consider the element
$$\sigma_{X,\Gamma_0(E)}(\widetilde{\Theta})\in \RKK^G(X; \Gamma_{\tau}(E),
\Gamma_0(E)\hot_{C_0(X)} C_{\tau}(U)[0,1]).$$ The balanced tensor
product $\Gamma_0(E)\hot_{C_0(X)} C_{\tau}(U)[0,1]$ is canonically
isomorphic to the algebra $C_\tau(U)\hot_{C_0(U)}\Gamma_0(P^*E)$,
where $P(x,y,t)= p_t(x,y)$ is the map of Notations \ref{P}.

Second, let $\varpi:\widetilde{U}\times[0,1]\longrightarrow
{U}\times[0,1];\,(x,y,t)\mapsto (y,x,1-t).$ Then this induces a
$G$-equivariant isomorphism
$$C_\tau(U)\hot_{C_0(U)}\Gamma_0(P^*E)\cong C_{\tau}(\widetilde{U})\hot_{C_0(\widetilde{U})}\Gamma_0(\varpi^*P^*E)$$
and observe
that $(P^*E)_{(x,y)}=E_{\gamma_{xy}(t)}$ and
$\left(\varpi^*P^*E\right)_{(x,y)}=E_{\gamma_{yx}(1-t)}$.
The isomorphism $\varphi: P^*E\to P_0^*E$ turns into an isomorphism
$\widetilde{\varphi}:\varpi^*P^*E\to \varpi^*P_0^*E=(P_0\circ \varpi)^*E$, and we notice that
$P_0\circ \varpi:\widetilde{U}\times [0,1]\to X$ equals $P_1:\widetilde{U}\times [0,1]\to X; P_1(x,y,t)=y$.
It follows that
\begin{align*}
C_\tau(U)\hot_{C_0(U)}\Gamma_0(P^*E)&\cong
C_{\tau}(\widetilde{U})\hot_{C_0(\widetilde{U})}\Gamma_0(\varpi^*P^*E)\\
&\cong
C_{\tau}(\widetilde{U})\hot_{C_0(\widetilde{U})}\Gamma_0(\varpi^*P_0^*E)\\
&\cong
C_{\tau}(\widetilde{U})\hot_{C_0(\widetilde{U})}\Gamma_0(P_1^*E)\\
&\cong
C_{\tau}(\widetilde{U})\hot_{C_0(\widetilde{U})}\Gamma_0(p_1^*E)[0,1]\\
&\stackrel{\id\hot \tilde\phi_1^{-1}}{\cong}
C_{\tau}(\widetilde{U})\hot_{C_0(\widetilde{U})}\Gamma_0(p_0^*E)[0,1],
\end{align*}
where $\tilde\phi_1:\Gamma_0(p_0^*E)\to \Gamma_0(p_1^*E)$ with respect to
the base $\widetilde{U}$
is induced  by the isomorphism $\phi_1:\Gamma_0(p_1^*E)\to \Gamma_0(p_0^*E)$ on the base $U$
via the flip $\sigma:U\to \widetilde{U}$.

Denote by $\Psi$ the above chain of isomorphisms. Note that $\Psi$ maps the
fibre $\Cl_{x,y}\hot E_x$ over the point $(x,y,0)\in U\times [0,1]$ identically to the fibre
$\Cl_{x,y}\hot E_x$ over $(y,x,1)\in \widetilde{U}\times[0,1]$ and the fibre
$\Cl_{x,y}\hot E_y$ over $(x,y,1)\in U\times[0,1]$ to the fibre
$\Cl_{x,y}\hot E_x$ over $(y,x,0)\in \widetilde{U}\times[0,1]$ via the map
$\id\hot \varphi_{y,x,1}^{-1}$.

Now, forgetting the $C_0(X)$-structure of $\sigma_{X,\Gamma_0(E)}(\widetilde{\Theta})$, we get
\begin{align*}
\Psi_*\left(\sigma_{X,\Gamma_0(E)}(\widetilde{\Theta})\right)
\in  \KK^G&\left(\Gamma_{\tau}(E),
C_{\tau}(\widetilde{U})\hot_{C_0(\widetilde{U})}\Gamma_0(p_0^*E)[0,1]\right)\\
&\quad\quad\quad \stackrel{\iota_*}{\to}
\KK^G\big(\Gamma_{\tau}(E),
\left(\Gamma_{\tau}(E)\hot C_\tau(X)\right)[0,1]\big).
\end{align*}
If we evaluate this class at $0$ and $1$ we get the triples:
$$ \left(C_\tau(\tilde{U})\hot_{C_0(X),p_0}\Gamma_0(E), \widetilde{\psi}_0,\widetilde{\Theta}_0\right)
\quad \text{and}\quad
 \left(C_\tau(\tilde{U})\hot_{C_0(X),p_0}\Gamma_0(E), \widetilde{\psi}_1,\widetilde{\Theta}_1\right)$$
 with
 $$\widetilde{\Theta}_0(x,y)= \frac{\rho(x,y)}{r(y)}\left(d_y\rho\right)(x,y)
 \quad\text{and}\quad
 \widetilde{\Theta}_1(x,y)= \frac{\rho(x,y)}{r(y)}\left(d_x\rho\right)(x,y)$$
 and the left actions $\widetilde{\psi}_0$ and $\widetilde{\psi}_1$ given by
 \begin{align*}
 \left(\widetilde{\psi}_0(f\hot_{C_0(X)}\eta)(g\hot \xi)\right)(x,y)
& =f(x)g(x,y)\hot \eta(x)\xi(x)\quad\text{and}\\
 \left(\widetilde{\psi}_1(f\hot_{C_0(X)}\eta)(g\hot \xi)\right)(x,y)
 &= f(y)g(x,y)\hot \varphi^{-1}_{y,x,1}(\eta(y))\xi(x),
 \end{align*}
for $f\hot_{C_0(X)}\eta\in
C_{\tau}(X)\hot_{C_0(X)}\Gamma_0(E),\,g\hot\xi\in
C_{\tau}(\widetilde{U})\hot_{C_0(X),p_0}\Gamma_0(E)$, where we have
carefully evaluated the isomorphism $\Psi$ on the fibres over
$x,y\in U$. Thus we see that evaluation at $1$ yields the element
$\Sigma_*\left(\sigma_{X,C_{\tau}(X)}(\Theta_E)\right)$ and
evaluation at $0$ gives an element which differs from $\sigma_{X,
\Gamma_0(E)}(\Theta_X)$ only by the fact that $U$ is replaced by
$\widetilde{U}$ and that
$\Theta(x,y)=\frac{\rho(x,y)}{r(x)}\left(d_y\rho\right)(x,y)$ is
replaced by $\widetilde{\Theta}_0$ as described above. If
$r(x)=\eta$ is constant (which we may always assume if $X$ is
$G$-compact), then $U=\widetilde{U}$ and $r(x)=r(y)$ for all
$(x,y)\in U$ and we are done. Otherwise, similar to the second part
of the proof of \cite[Lemma 4.6]{Kas1} we deform $U$ into
$\widetilde{U}$ and $\Theta$ into $\widetilde{\Theta}_0$ to finish
the proof.
\end{proof}

We are finally ready for

\begin{proof}[Proof of Theorem \ref{thm-kas-bundle}]
We have to show that the maps $\mu$ and $\nu$ of Theorem
\ref{thm-kas-bundle} are inverse to each other.
\begin{eqnarray*}
\lefteqn{(\nu\circ\mu)(\alpha)=[\Theta_E]\hot_{C_0(X)\hot\Gamma_{\tau}(E)}\sigma_{C_0(X)}\left(\sigma_{X,C_{\tau}(X)}(\alpha)\hot_{C_{\tau}(X)}[d_X]\right)}\\
&=&\left([\Theta_E]\hot_{C_0(X)\hot\Gamma_{\tau}(E)}\sigma_{C_0(X)}\sigma_{X,C_{\tau}(X)}(\alpha)\right)\hot_{C_0(X)\hot
C_{\tau}(X)}\sigma_{C_0(X)}[d_X]\\
&=&\left(\alpha\hot_{C_0(X)}[\Theta_X]\right)\hot_{C_0(X)\hot
C_{\tau}(X)}\sigma_{C_0(X)}[d_X]\hskip2cm \text{by Lemma
\ref{lem4.5}}\\
&=&\alpha\hot_{C_0(X)}\left([\Theta_X]\hot_{C_0(X)\hot
C_{\tau}(X)}\sigma_{C_0(X)}[d_X]\right)\\
&=&\alpha \hskip2cm \text{by Lemma 4.8 in \cite{Kas1}}.
\end{eqnarray*}

\begin{eqnarray*}
\lefteqn{(\mu\circ\nu)(\beta)=\sigma_{X,C_{\tau}(X)}\left([\Theta_E]\hot_{C_0(X)\hot\Gamma_{\tau}(E)}\sigma_{C_0(X)}(\beta)\right)\hot_{C_{\tau}(X)}[d_X]}\\
&=&\sigma_{X,C_{\tau}(X)}[\Theta_E]\hot_{C_{\tau}\hot
\Gamma_{\tau}(E)}\left(\sigma_{C_{\tau}(X)}(\beta)\hot_{C_{\tau}(X)}[d_X]\right)\\
&=&\sigma_{X,C_{\tau}(X)}[\Theta_E]\hot_{C_{\tau}\hot
\Gamma_{\tau}(E)}\left(\beta\hot[d_X]\right)\\
&=&\sigma_{X,C_{\tau}(X)}[\Theta_E]\hot_{C_{\tau}\hot
\Gamma_{\tau}(E)}\left([d_X]\hot\beta\right)\\
&=&\sigma_{X,C_{\tau}(X)}[\Theta_E]\hot_{C_{\tau}\hot
\Gamma_{\tau}(E)}\left(\sigma_{\Gamma_{\tau}(E)}[d_X]\hot_{\Gamma_{\tau}(E)}\beta\right)\\
&=&\left(\sigma_{X,\Gamma_{\tau}(E)}[\Theta_X]\hot_{
\Gamma_{\tau}(E)\hot C_{\tau}(X)}\sigma_{\Gamma_{\tau}(E)}[d_X]\right)\hot_{\Gamma_{\tau}(E)}\beta\hskip1cm\text{by Lemma \ref{lem4.6}}\\
&=&\sigma_{X,\Gamma_{\tau}(E)}\left([\Theta_X]\hot_{C_0(X)\hot
C_{\tau}(X)}\sigma_{C_0(X)}[d_X]\right)\hot_{\Gamma_{\tau}(E)}\beta\\
&=&\beta\hskip2cm\text{by Lemma 4.8 in \cite{Kas1}}.
\end{eqnarray*}
\end{proof}

\section{Poincar\'e-duality for crossed products by proper actions}\label{sec-prop}

Suppose that $G$ is a locally compact group which acts on the locally compact space
$X$.   The $G$-equivariant Brauer group $\Br_G(X)$ of $X$ in the sense of \cite{CKRW}
 consists of $X\rtimes G$-equivariant Morita
equivalence classes of continuous-trace
C*-algebras $A$ with base $X$ and which are equipped with an action of  $G$ which covers the
given action on $X$.
 The group operation in $\Br_G(X)$  is
 given by taking tensor products over $X$ and diagonal actions.
We shall write
 $\delta$ for a class in $\Br_G(X)$, and we shall denote by $C_0(X,\delta)$
 a representative of the corresponding continuous trace algebra equipped with the appropriate action of $G$. Moreover,
 we write $C_0(X, \delta^{-1})$ for a representative of the inverse $\delta^{-1}\in \Br_G(X)$,
 and we write $C_0(X, \delta\cdot\mu)$ for a representative of the product $\delta\cdot \mu\in \Br_G(X)$
 if $\delta,\mu\in \Br_G(X)$. In this notation we have (up to $X\rtimes G$-equivariant Morita equivalence):
 \begin{equation}\label{eq-inverse}
 C_0(X,\delta)\hot_{C_0(X)}C_0(X,\mu)=C_0(X,\delta\cdot \mu)
 \quad\text{and}\quad
 C_0(X,\delta)\hot C_0(X,\delta^{-1})=C_0(X),
 \end{equation}
as $X\rtimes G$-algebras.

\begin{lem}\label{lem-inverse}
If $\delta\in \Br_G(X)$, then for every pair of $X\rtimes G$-algebras $A$ and $B$,
 tensoring over $C_0(X)$ with $C_0(X,\delta)$ gives an
isomorphism
$$\sigma_{X,C_0(X,\delta)}:\RKK^G(X; A, B)\stackrel{\cong}{\longrightarrow}
\RKK^G\big(X; C_0(X,\delta)\hot_{C_0(X)}A, C_0(X,\delta)\hot_{C_0(X)}B\big).$$
In particular, there is a canonical isomorphism
$$\RKK^G\big(X; C_0(X,\delta)\hot_{C_0(X)}A, B\big)
\cong
\RKK^G\big(X; A, C_0(X,\delta^{-1})\hot_{C_0(X)}B\big).$$
\end{lem}
\begin{proof} By (\ref{eq-inverse}) it is instantly clear that $\sigma_{X,C_0(X,\delta^{-1})}$ is
an inverse to $\sigma_{X,C_0(X,\delta)}$. The second isomorphism follows easily from the first.
\end{proof}

If $X$ is a $G$-manifold, then we shall say that a class $\delta\in
\Br_G(X)$ is {\em smooth}, if it can be represented by a smooth
$G$-equivariant locally trivial bundle $p:E\to X$ of elementary
C*-algebras equipped with a smooth action of $G$.
Moreover, throughout this section we shall use the notation
$C_0(X,\delta\cdot\tau):=C_0(X,\delta)\hot_{C_0(X)}\Ctau$.

As a direct consequence of Theorem \ref{thm-kas-bundle}
and the above lemma we obtain
an equivariant version of Poincar\'e duality  in twisted $K$-theory, i.e., for
$G$-equivariant continuous-trace algebras over $X$.

\begin{cor}\label{cor-compact}
Assume that $G$ is a compact group, $X$ is a compact  $G$-manifold and
$\delta\in \Br_G(X)$ is smooth.
Then, for all $G$-algebras $A$ and $B$ there are natural (in $A$ and $B$) isomorphisms
$$\Phi_{A,B}:\KK^G(C_0(X,\delta\cdot\tau)\hot A, B)\stackrel{\cong}{\longrightarrow}
\KK^G(A, C_0(X,\delta^{-1})\hot B).$$
\end{cor}
\begin{proof}
Note that with our notations, we have $\Gamma_0(E)=C_0(X,\delta)$ and
$\Gamma_0(E)\hot_{C_0(X)}\Ctau=C_0(X,\delta\cdot\tau)$
if $p:E\to X$ is a smooth bundle representing $\delta$.
By Proposition \ref{lem-trans} we know that $p:E\to X$ is feasible.
Thus, by Theorem \ref{thm-kas-bundle} and Lemma \ref{lem-inverse}
we have natural  isomorphisms
\begin{align*}
\KK^G(C_0(X, \delta\cdot\tau)\hot A,B)&\cong \RKK^G(X; C_0(X,\delta)\hot A, C_0(X,B))\\
&\cong \RKK^G(X; C_0(X, A), C_0(X,\delta^{-1})\hot B)\\
&= \KK^G(A, C_0(X,\delta^{-1})\hot B),
\end{align*}
where the last identification can be made because $X$ is compact.
\end{proof}

\begin{remark}\label{rem-Tu} As mentioned in the introduction,
basically the same result (for compact Lie groups but 
without the smoothness assumption
given above) was stated independently by Jean-Louis Tu in his
recent preprint \cite{Tu-P}. Although he does not make this explicit, his 
argument uses a certain locality condition for the bundle as spelled 
out in \cite[\S 6]{AS}. 
\end{remark}

We next recall the notion of Poincar\'e duality for $C^*$-algebras,
see for example \cite{Em}.

\begin{defn} \rm Let $\Lambda$ and $\widehat{\Lambda}$ be graded
$C^*$-algebras. Then $\Lambda$ and $\widehat{\Lambda}$ are
\emph{Poincar\'e dual} if there exist classes $\Delta \in
\KK(\Lambda \hot \widehat{\Lambda} , \C)$ and $\widehat{\Delta} \in
\KK(\C, \widehat{\Lambda}\hot\Lambda )$ such that
\begin{equation}\label{fundamentalclass}
\widehat{\Delta}\hot_{\Lambda} \Delta = 1_{\widehat{\Lambda}},
\;\;\;\;\;\; \widehat{\Delta}\hot_{\widehat{\Lambda}} \Delta =
1_{\Lambda}. \end{equation}
Equivalently, we are given a system of
isomorphisms
$$\KK(\Lambda \hot A, B) \stackrel{\cong}{\longrightarrow}
\KK(A, \widehat{\Lambda} \hot B)$$
natural with respect to intersection and composition products.
\end{defn}

\begin{remark}\label{symmetry}
Note that when we say $\Lambda$ and $\widehat{\Lambda}$ are
\emph{Poincar\'e dual}, we already implicitly used the fact that
Poincar\'e duality is symmetric. Indeed one can show that
$\Delta':=\Sigma^*(\Delta)\in \KK(\widehat{\Lambda}\hot\Lambda,\C)$
and $\widehat{\Delta}':=\Sigma_*(\widehat{\Delta})\in
\KK(\C,\Lambda\hot\widehat{\Lambda})$ satisfy Equation
(\ref{fundamentalclass}), where
$\Sigma:\widehat{\Lambda}\hot\Lambda\longrightarrow\Lambda\hot\widehat{\Lambda}$
is the flip isomorphism. Likewise the equivariant Poincar\'e duality
(i.e., the above definition with the functor $\KK$ replaced by
$\KK^G$) is also symmetric.
\end{remark}

We now consider the case where a {\em discrete} group $G$ acts
 {\em properly} on the $G$-manifold $X$, and that, in addition, $G\backslash X$
 is compact. Then, if $\delta\in \Br_G(X)$ is smooth, we will show that the
 algebras $C_0(X,\delta\cdot\tau)\rtimes G$ and $C_0(X,\delta^{-1})\rtimes G$
 are Poincar\'e dual. For this we need some additional tools.
 The following lemma is a part  of the proof of \cite[Theorem 5.4]{KasSk}
 (a generalization to
a much weaker notion of proper actions is given by Ralf Meyer in \cite{M}):

\begin{lem}\label{lem-propermodule}
Suppose that $X$ is a proper $G$-space and that $D$ is a $X\rtimes G$-algebra.
Then there is a natural equivalence between the category
of $G$-equivariant Hilbert $D$-modules and the category of Hilbert $D\rtimes G$-modules
sending a Hilbert  $D$-module $\E$ to the Hilbert $D\rtimes G$-module $\widetilde{\E}$
constructed as follows:
Write $\mathcal E_c:= C_c(X)\cdot \mathcal E$ ($:=\E\cdot C_c(X)$ with right action of $C_c(X)$
on $\E$ given via the right action of $D$ and the structure map $C_0(X)\to ZM(D)$).
Then for $\xi, \eta \in \mathcal E_c$ and
$f\in C_c(G,D)\subseteq D\rtimes G$ define a  $C_c(G,D)$-valued inner product
and a right module action of $C_c(G,D)$ on $\E_c$ by:
\begin{equation}\label{eq-module}
\begin{split}
\lk \xi,\eta\rk_{C_c(G,D)}(s)&=\big(s\mapsto \sqrt{\Delta_G(s^{-1})}\lk \xi,\gamma_s( \eta)\rk_D\big)
\quad\quad\text{and}\\
\xi\cdot f&=\int_G\sqrt{\Delta_G(s^{-1})} \gamma_s(\xi\cdot f(s^{-1}))\,ds,
\end{split}
\end{equation}
where $\gamma:G\to\Aut(\E)$ denotes the given action on $\E$.
Then  $\mathcal E_c$ completes to the corresponding
Hilbert $D\rtimes G$-module $\widetilde{\E}$.
The algebra of adjointable operators $\mathcal L(\widetilde{\E})$
is then isomorphic to the algebra $\mathcal L(\E)^G$ of $G$-invariant operators in
$\mathcal L(\mathcal E)$ and $\mathcal K(\widetilde E)$
coincides with the generalized fixed-point algebra $\mathcal K(\mathcal E)^G$
in the sense of \cite[Section 3]{Kas1}.
\end{lem}

\begin{remark}\label{rem-Mischenko}
In the special case where $D=C_0(X)$ and where $\E=C_0(X)$ is viewed
as a $C_0(X)$-module in the canonical way, the above lemma gives a
corresponding Hilbert $C_0(X)\rtimes G$-module
$\E_X:=\widetilde{\E}=\overline{C_c(X)}$. In this case
$C_0(G\backslash X)$ is the generalized fixed-point algebra of
$C_0(X)=\K(\E)$, and hence we have $C_0(G\backslash X)=\K(\E_X)$.
The pair $(\E_X, 0)$ determines  a canonical class
$\Lambda_X\in\RKK(G\backslash X; C_0(G\backslash X), C_0(X)\rtimes
G)$ and, if $G\backslash X$ is compact, a corresponding class
$\lambda_X\in \KK(\CC, C_0(X)\rtimes G)$. The class $\Lambda_X$
(resp. $\lambda_X$) is often called the
 {\em Mishchenko line bundle for $X$}.
\end{remark}

The next  result is a slight extension of a result of Kasparov and Skandalis
 (\cite[Theorem 5.4]{KasSk}) which treats the case where $A=\CC$.
 The case where $A$ is arbitrary and $D=C_0(X,B)$
for some $G$-algebra $B$ was treated by
Emerson and  Meyer (see \cite[Lemma 20]{EM}). Since the proof is a direct extension
of the proofs of those special cases, we restrict ourselves to explain the basic
steps:

 \begin{thm}[{cf. \cite[Theorem 5.4]{KasSk} and \cite[Lemma 20]{EM}}]\label{thm-EM}
 Suppose that $G$ is a locally compact group and $X$ is a proper $G$-compact $G$-space.
 Suppose further that $A$ is a $C^*$-algebra equipped with the trivial $G$-action and that $D$ is
an $X\rtimes G$-algebra. Then there is a natural isomorphism
 $$\Theta_{A,D}:\RKK^G(X; C_0(X,A),D)\stackrel{\cong}{\longrightarrow} \KK(A, D\rtimes G).$$
 given by the composition
 $$\RKK^G(X; C_0(X,A),D)\stackrel{J_G}{\longrightarrow} \KK\big((C_0(X)\rtimes G)\hot A, D\rtimes G\big)
 \stackrel{\sigma_A(\lambda_X)\hot}{\longrightarrow} \KK(A, D\rtimes G),$$
 where $J_G$ denotes Kasparov descent morphism (and using the fact that
 $C_0(X,A)\rtimes G\cong (C_0(X)\rtimes G)\hot A$, since $G$ acts trivially on $A$),
 and $\lambda_X\in K_0(C_0(X)\rtimes G)$ is the  Mishchenko line bundle for $X$.
 \end{thm}
\begin{proof}
Let  $(\mathcal E, T)$ be a cycle for  $\RKK^G(X; C_0(X,A), D)$.
Consider the operator $\widetilde{T}=\int_G s(cT)\,ds$,
where $c:X\to [0,\infty)$ is  any
compactly supported continuous function on $X$ such that $\int_G c(s^{-1}x)\,ds=1$ for all $x\in X$
(the arguments given
in the proof of \cite[Theorem 3.4]{Kas1} show that the $t\mapsto T+t(\tilde{T}-T)$ gives
 an operator homotopy between $T$ and $\widetilde{T}$, so that $[\E,\pi,T)]=[(\E,\pi,\widetilde{T})]$).
There is a canonical pairing $C_c(X,A)\times  \mathcal E\to \mathcal E_c$ which integrates to a
map $C_c(X,A)\odot C_c(G,\mathcal E) \to \mathcal E_c$ given by
$$\xi\hot g\mapsto F_{\xi,g}:=\int_G\sqrt{\Delta_G(s^{-1})} s^{-1}\cdot (\xi  \cdot (g(s)))\,ds$$
A straightforward computation shows that this map
preserves the inner products and the right actions of $C_c(G,D)$, so that it extends to an isomorphism
$$(\E_X\hot A)\hot_{(C_0(X)\rtimes G)\hot A} \E\rtimes G\stackrel{\cong}{\longrightarrow}
\widetilde{\E},$$
with $\widetilde{\E}$ as in Lemma \ref{lem-propermodule}.
We may then proceed precisely as in the proof of \cite[Theorem 5.4]{KasSk} to see that
$[(\E,T)]\to [(\widetilde{\E},\widetilde{T})]$ defines an isomorphism from
$\RKK^G(X; C_0(X,A),D)$ to $\KK(A, D\rtimes G)$, which coincides with
the map as given in the theorem.
\end{proof}

Combining this theorem with Theorem \ref{thm-kas-bundle} we are now able  to prove:

\begin{thm}\label{thm-poincare-bundle}
Suppose that $G$ is a discrete group acting properly on the
$G$-manifold $X$ such that $G\backslash X$ is compact and assume that $\delta\in \Br_G(X)$
is smooth.
Then
$C_0(X,\delta\cdot\tau)\rtimes G$ and $C_0(X,\delta^{-1})\rtimes G$
are Poincar\'e dual.

If, in addition, $X$ is equipped with  a $G$-equivariant spin$^c$ structure, then
there is a Poincar\'e duality between
$C_0(X,\delta)\rtimes G$ and  $C_0(X,\delta^{-1})\rtimes G$ of degree
$\dim(X)$ mod $2$, i.e.,
up to
a dimension shift of order $\dim(X)$ mod $2$.
\end{thm}

We show this using the natural isomorphism version of the definition: that
is, by showing existence of a natural system of isomorphisms
$$\Phi_{A,B}: \KK\big((C_0(X,\delta\cdot\tau)\rtimes G)\hot A, B)\stackrel{\cong}{\longrightarrow}
\KK(A, (C_0(X,\delta^{-1})\rtimes G)\hot B),$$
one for each $A,B$.

\begin{proof}
We define the map $\Phi_{A,B}$ as a composition
\begin{equation}\label{eq-comp}
\Phi_{A,B}:=E_{A,B}\circ F_{A,B}\circ G_{A,B}\circ H_{A,B},
\end{equation}
where
$$H_{A,B}:\KK\big((C_0(X,\delta\cdot\tau)\rtimes G)\hot A, B\big)
\stackrel{\cong}{\to} \KK^G(C_0(X,\delta\cdot\tau)\hot A, B)$$
is the canonical isomorphism due to discreteness of $G$.
The map
$$G_{A,B}: \KK^G(C_0(X,\delta\cdot\tau)\hot A, B)\stackrel{\cong}{\to}
\RKK^G(X; C_0(X,\delta)\hot A, C_0(X,B))$$
is the inverse of the isomorphism of Theorem \ref{thm-kas-bundle}.
The third map
$$F_{A,B}:\RKK^G(X; C_0(X,\delta)\hot A, C_0(X,B))\to
\RKK^G(X; C_0(X,A), C_0(X,\delta^{-1})\hot B)$$
is the isomorphism of Lemma \ref{lem-inverse}. Finally, the map
$$E_{A,B}:  \RKK^G(X; C_0(X,A), C_0(X,\delta^{-1})\hot B)\to
\KK\big(A, (C_0(X,\delta^{-1})\rtimes G)
\hot B\big)$$
is the map of Theorem \ref{thm-EM} combined with the canonical isomorphism
$$(C_0(X,\delta^{-1})\hot B)\rtimes G\cong  (C_0(X,\delta^{-1})\rtimes G)\hot B$$
which results from the fact that $G$ acts trivially on $B$.
It follows from the general properties of Kasparov products as
outlined in \cite{Kas1} that all maps above are natural in $A$ and
$B$ with respect to taking Kasparov products with elements in
$\KK(A',A)$ or $\KK(B, B')$, respectively. They are also natural
with respect of taking external Kasparov products (i.e.~Kasparov
products over $\CC$ with elements in $\KK(A', B')$).

The final assertion now follows from the fact that existence of a $G$-equivariant spin$^c$
structure implies the existence of a $X\rtimes G$-equivariant Morita equivalence
$$\Ctau\sim_M C_0(X)  \quad (\text{resp. $\Ctau\hot Cl_1\sim_M C_0(X)$}),$$
if $X$ has even (resp. odd) dimension (where $Cl_1$ denotes the first Clifford algebra
equipped with the trivial $G$-action).
\end{proof}

\begin{remark}\label{rem-bundle} The above theorem can easily be extended to the
following more general situation: Assume as in the theorem that $G$ is discrete
and that $X$ is a proper $G$-manifold such that $G\backslash X$ is compact.
Assume further that $p:E\to X$ is a feasible locally trivial bundle of C*-algebras
over $X$ (where we do {\bf not} assume that the fibres are elementary C*-algebras)
such that there exists a $C_0(X)$-algebra $\Gamma_0(\bar{E})$ (this is just convenient
notation---we do not assume that $\bar{E}$ is a locally trivial or even continuous bundle of
C*-algebras over $X$) which is inverse to $\Gamma_0(E)$ in $\RKK^G(X;\cdot,\cdot)$, i.e.,
we have
$$\Gamma_0(E)\hot_{C_0(X)}\Gamma_0(\bar{E})\sim_{\RKK^G}C_0(X).$$
Then replacing $C_0(X,\delta\cdot \tau)$ by
$\Gamma_{\tau}(E):=\Gamma_0(E)\hot_{C_0(X)}\Ctau$ and $C_0(X,\delta^{-1})$
by $\Gamma_0(\bar{E})$ in the proof of the above theorem will show that
$\Gamma_{\tau}(E)\rtimes G$ is Poincar\'e dual to $\Gamma_0(\bar{E})\rtimes G$.
\end{remark}

\section{Some Applications}\label{sec-twistedgroup}
Throughout this section we want to assume that $G$ is a discrete group such
that there exists
a $G$-manifold $X$  which satisfies the following axioms:
\begin{itemize}
\item[(A1)] $G$ acts properly on $X$ and  $G\backslash X$ is compact.
\item [(A2)] $X$ is a special $G$-manifold as defined in \cite[Definition 5.1]{Kas1}.
\item[(A3)] $X$ is $G$-equivariantly spin$^c$.
\end{itemize}
Notice that (A2) implies the existence of an element $\eta_X\in \KK^G(\CC,C_\tau(X))$,
unique under the conditions spelled out in \cite[Definition 5.1]{Kas1},
which is dual to the Dirac element $d_X\in\ \KK^G(C_{\tau}(X),\CC)$ in the sense that
\begin{equation}\label{eq-dual-Dirac}
d_X\hot\eta_X=1_{C_\tau(X)}\in \KK^G(C_{\tau}(X),C_{\tau}(X)).
\end{equation}
The reverse product
$$\gamma_G=\eta_X\hot_{C_{\tau}(X)}d_X\in \KK^G(\CC,\CC)$$
is then independent of the special choice of $X$, and it is called the {\em $\gamma$-element}
of $G$.
We extend the above list of axioms by
\begin{itemize}
\item[(A4)] $\gamma_G=1\in \KK^G(\CC,\CC)$
\end{itemize}

Any cocompact discrete subgroup $G$ of an (almost)
connected group $L$ satisfies axioms (A1) and (A2) with
$X=L/K$, the quotient of $L$ by the maximal compact
subgroup $K$ equipped with the canonical
$L$-equivariant metric.
Notice that in this case, modulo passing to a double cover $\widetilde{L}$ of $L$,
which, unfortunately, we have to pay by passing also to a double cover $\widetilde{G}$ of our
subgroup $G\subseteq L$, the $\spin^c$ axiom (A3) can always be arranged.
However, Axiom (A4) imposes a more severe restriction on the groups, but it applies
to all amenable (or, more generally, a-$T$-menable) groups $G$ (see \cite{HK}).

If $\gamma_G=1$, it follows that $C_\tau(X)$ is $\KK^G$-equivalent
to $\CC$, which then implies via Kasparov's descent homomorphism that
\begin{equation}\label{eq-KKequiv}
(A\hot C_{\tau}(X))\rtimes_{(r)}G\sim_{\KK} A\rtimes_{(r)}G
\end{equation}
for the full and reduced crossed products of $A\hot
C_{\tau}(X)\sim_{\KK^G}A$ by $G$ and for any $G$-C*-algebra $A$.
Since $G$ acts properly on $X$, we have $(A\hot C_{\tau}(X))\rtimes
G\cong (A\hot C_{\tau}(X))\rtimes_{r}G$ by \cite[Theorem
3.13]{Kas1}, and hence we see that $G$ is $\KK$-amenable in the
sense that $A\rtimes G\sim_{\KK} A\rtimes_r G$ via the quotient map
$A\rtimes G\to A\rtimes_rG$.

We use these observations to prove

\begin{thm}\label{thm-spin}
Suppose that $G$ satisfies Axioms (A1), (A2), (A3) and (A4). Assume further that $M$ is a compact
$G$-manifold  and that $\delta\in \Br_G(M)$ is smooth. Then
$C(M,\delta\cdot\tau)\rtimes G$ and $C(M,\delta^{-1})\rtimes G$ are Poincar\'e dual
of degree $\dim(X) \mod 2$.

If, in addition, $M$ has a $G$-equivariant spin$^c$-structure, then
$C(M,\delta)\rtimes G$ is  Poincar\'e dual to $C(M,\delta^{-1})\rtimes G$
of degree $\dim(X)+\dim(M) \mod 2$.
\end{thm}

\begin{remark} In Theorem
\ref{thm-poinc-comp} below we shall relax
considerably the condition (A4) used  above.
However, in that theorem we have to give an extra spin$^c$ assumption
which we can avoid in Theorem \ref{thm-spin} above.
\end{remark}

\begin{proof}[Proof of Theorem \ref{thm-spin}]
By Axiom (A3) we know that $\Ctau$ (resp. $\Ctau\hot Cl_1$ if
$\dim(X)$ is odd) is $X\rtimes G$-equivariantly Morita equivalent to $C_0(X)$.
We therefore obtain a $\KK^G$-equivalence of degree $\dim(X) \mod 2$
as in (\ref{eq-KKequiv}) with $\Ctau$ replaced by $C_0(X)$.
Moreover,  by the general properties of Clifford-bundles
we have $\ctau(X\times M)\cong \Ctau\hot \ctau(M)$.
Assume now that $\delta\in \Br_G(M)$ is smooth and let $p:X\times M\to M$ denote the
canonical projection. Then $C_0(X\times M,p^*\delta)=C_0(X)\hot C(M,\delta)$
and $C_0(X\times M, p^*\delta\cdot\tau)=\Ctau\hot C(M,\delta\cdot\tau)$,
where we use the notation of Theorem \ref{thm-poincare-bundle}.
Thus, (\ref{eq-KKequiv}) and its variant for $C_0(X)$ imply that
$$C_0(X\times M, p^*\delta\cdot\tau)\rtimes G\sim_{\KK} C(M,\delta\cdot\tau)\rtimes G
\quad\text{and}\quad
C_0(X\times M, p^*\delta^{-1})\rtimes G\sim_{\KK} C(M,\delta^{-1})\rtimes G,$$
where the latter is of degree $\dim(X)\mod 2$.
Theorem \ref{thm-poincare-bundle} applied to the smooth class $p^*\delta\in \Br_G(X\times M)$
with inverse $(p^*\delta)^{-1}=p^*(\delta^{-1})$ provides a Poincar\'e duality
between
\mbox{$C_0(X\times M, p^*\delta\cdot\tau)\rtimes G$} and
$C_0(X\times M, p^*\delta^{-1})\rtimes G$ which, by $\KK$-equivalence,
implies duality of degree $\dim(X)\mod 2$ between
$C(M,\delta\cdot\tau)\rtimes G$ and $C(M,\delta^{-1})\rtimes G$.

The statement for the case where $M$ is also spin$^c$ is now obvious.
\end{proof}

In what follows next, we want to specialize the above result to the case where $M=\{\pt\}$.
In that case the algebra $C(M,\delta)$ will just be the algebra of compact operators
on some Hilbert space equipped with an action of $G$. Such actions are classified
up to Morita equivalence by the second group cohomology $H^2(G,\TT)$ with coefficient
the circle group $\TT$ considered as trivial $G$-module.
To be more precise, if $\K=\K(H)$ and $\alpha:G\to \Aut(\K)$ is given, one can choose a map
 $V:G\to \U(H)$ such that $\alpha_s=\Ad V_s$ for all $s\in G$. Then there is a cocycle
$\om_{\alpha}\in Z^2(G,\TT)$ determined by the equation
\begin{equation}\label{eq-V}
V_s V_t=\om_{\alpha}(s,t)V_{st}\quad\text{for all $s,t\in G$}
\end{equation}
and $[\alpha]\mapsto [\om_\alpha]; \Br_G(\pt)\to H^2(G,\TT)$ is an isomorphism of groups
(e.g. see \cite{CKRW}).

On the other hand, given $\om \in Z^2(G,\TT)$, then we can construct
full and reduced  twisted group algebras $C^*(G,\om)$ and $C_r^*(G,\om)$ as follows:
The twisted convolution algebra
$\ell^1(G,\omega)$ is defined as the vector space of all summable
complex functions on $G$ with convolution and involution given by
\begin{eqnarray*}f\ast_{\omega}g(s)&:=&\sum_{t\in G} f(t)g(t^{-1}s)\omega(t,t^{-1}s)\\
f^*(s)&:=&\overline{\omega(s,s^{-1})f(s^{-1})}.
\end{eqnarray*}
If $V:G\to U(H)$ is any $\om$-representation of $G$, i.e., $V$ satisfies equation (\ref{eq-V}),
then $V$ determines a $\ast$-homomorphism $\tilde{V}$ of
$\ell^1(G,\omega)$ into $B(H)$ via the formula
$$\tilde{V}(f):=\sum_{s\in G}f(s)V(s), \quad f\in \ell^1(G,\omega),$$
and every nondegenerate $*$-representation of $\ell^1(G,\omega)$ appears in
this way. The full twisted group algebra $C^*(G,\omega)$ is defined
as the enveloping $C^*$-algebra of $\ell^1(G,\omega)$ and the reduced
twisted group algebra $C^*_r(G,\omega)$ is obtained from the regular
$\omega$-representation $L_{\omega}:G\longrightarrow U(\ell^2(G))$
given by
$$\left(L_{\omega}(s)\xi\right)(t):=\omega(s,s^{-1}t)\xi(s^{-1}t),~~~\xi\in
\ell^2(G),~s,t\in G.$$
Notice that the isomorphism classes of $C_{(r)}^*(G,\om)$ only depend on the class
$[\om]\in H^2(G,\TT)$!

Assume now that  $V$ is any $\bar{\om}$-representation on some Hilbert space $H$,
where $\bar{\om}$ denotes the {\em inverse} of the cocycle $\om$ and let
$\alpha=\Ad V:G\to \Aut(\K(H))$.
 Then  there is a canonical isomorphism
\begin{equation}\label{eq-isokomp}
\K \hot C_{(r)}^*(G,{\om})\cong \K\rtimes_{\alpha, (r)}G\quad \text{with $(\K=\K(H))$}
\end{equation}
for the full and reduced crossed products, given on the level of $\ell^1$-algebras by the formula
$k\hot f\mapsto (s\mapsto f(s)kV_s^*)$
(this  can be regarded as  a very special case of the stabilization theorem of \cite{PR}).
Using all this, we get as a special case of Theorem \ref{thm-spin}

\begin{thm}\label{thm-twist}
Suppose that $G$ satisfies the axioms (A1), (A2), (A3) and (A4).
 Then, for each $\om\in Z^2(G,\TT)$,
 $C^*(G,\bar{\om})$ and $C^*(G,\om)$  are Poincar\'e dual
of degree $j=\dim(X) \mod 2$.
\end{thm}
\begin{proof}
 Theorem \ref{thm-spin} applied to $M=\{\pt\}$ provides a Poincar\'e duality
for  $\K\rtimes_{\alpha_\om}G$ and $\K \rtimes_{\alpha_{\om^{-1}} } G$, where
we denote by $\alpha_{\om}$ an action corresponding to the class $[\om]\in H^2(G,\TT)$
as explained above. The result then follows from (\ref{eq-isokomp}).
  \end{proof}

Note that in  case $G=\ZZ^n$ the above theorem implies the well known self dualities
for the non-commutative $n$-tori, which are just
the twisted group algebras $C^*(\ZZ^n,\om)$ (in this case there is a canonical
isomorphism $C^*(\ZZ^n,\om)\cong C^*(\ZZ^n,\bar\om)$).

\begin{remark} Since the assumption made on $G$ in Theorems \ref{thm-spin} and
\ref{thm-twist} imply $K$-amenability of $G$, it is clear that the full crossed products
(resp. twisted group algebras)
in the statements of these theorems can be replaced
by the appropriate reduced crossed products (resp. twisted group algebras).
\end{remark}

We want to finish this section by a discussion on how one could
relax the assumptions in Theorem \ref{thm-spin}. In particular we
would like to get rid of the relatively strong assumption (A4). So
let us assume that $G$ satisfies axioms (A1), (A2), and (A3). The
idea is to replace  Axiom (A4) by
 the considerably weaker assumption that the $\gamma$-element
$\gamma_{M\rtimes G}$ of the groupoid $M\rtimes G$ is equal to
$1_{C(M)}\in \RKK^G(M;
C(M),C(M))$ ($ =\KK^{M\rtimes G}(C(M), C(M))$)
in the sense of Tu (see \cite[Section 5.2]{Tu1}).

Note that by the main result of \cite{Tu2}, this is always true if $G$ acts amenably
on $M$. But in order to get really more general results, we should {\bf not} assume
that $M$ admits a $G$-invariant metric, since this together with amenability
of the action of $G$ on $M$  would imply amenability
of $G$.

By the properties of the $\gamma$-element, as formulated in
\cite[Proposition 5.20]{Tu1}, one easily checks that
$\gamma_{M\rtimes G}=\sigma_{C(M)}(\gamma_G)$, and hence the Dirac
and dual-Dirac elements for $M\rtimes G$ are given by $d_{M\rtimes
G}=\sigma_{C(M)}(d_X)$ and $\eta_{M\rtimes
G}:=\sigma_{C(M)}(\eta_X)$, with $d_X$, $\eta_X$ as in the previous
discussions. By assumption we have
$$\eta_{M\rtimes G}\otimes d_{M\rtimes G}=\gamma_{M\rtimes G}=1_{C(M)}\in \RKK^G(M; C(M), C(M)).
$$
Hence, if $\delta\in \Br_G(M)$ is any element in the equivariant
Brauer group of $M$, then tensoring everything above with
$C(M,\delta)$ over $C(M)$ provides a $\RKK^G$-, and hence also a $
\KK^G$-equivalence $\tilde{D}:=\sigma_{M, C(M,\delta)}(d_{M\rtimes
G})$ between $C_0(X\times M, p^*(\delta))$ and $C(M,\delta)$, and
hence $KK$-equivalence between the crossed products $C_0(X\times M,
p^*(\delta))\rtimes G$ and $C(M,\delta)\rtimes G$, where $p:X\times
M\to M$ denotes the projection. Note that, since full and reduced
crossed products coincide for proper actions, we may replace
$C(M,\delta)\rtimes G$ by the reduced crossed product
$C(M,\delta)\rtimes_rG$ if we wish.

Since $G$  acts properly on $X\times M$ we may choose a
$G$-invariant metric on $X\times M$. We assume that this metric
admits a $G$-invariant spin$^c$-structure and also assume that
$\delta\in \Br_{G}(M)$ is smooth. Then we may apply Theorem
\ref{thm-poincare-bundle} to obtain a $\KK$-theoretic Poincar\'e
duality between $C_0(X\times M, p^*(\delta))\rtimes G$ and
$C_0(X\times M, p^*(\delta)^{-1})\rtimes G$, which by the above
observed $\KK$-equivalences provides Poincar\'e dualities for
$C(M,\delta)\rtimes_{(r)}G$ and $C(M,\delta^{-1})\rtimes_{(r)}G$
where the subscript $(r)$ indicates that we may take maximal or
reduced crossed products at any side as we wish!

So, putting things together, we obtain

\begin{thm}\label{thm-poinc-comp} Suppose that $G$ satisfies axioms (A1), (A2), (A3)
 above with respect to the proper $G$-manifold $X$.
Assume that $M$ is  any compact $G$-manifold  (we {\bf do not}
assume that $G$ acts isometrically on $M$) such that the
$\gamma$-element of the groupoid $M\rtimes G$ in the sense of Tu
(see \cite{Tu1}) is equal to $1_{C(M)}$ (which is automatic if $G$
acts amenably on $M$ by \cite{Tu2}). Assume further that $X\times M$
admits a $G$-equivariant spin$^c$-structure.  Then, for any smooth
class $\delta\in \Br_G(M)$, there exists a $\KK$-theoretic
Poincar\'e duality for $C(M,\delta)\rtimes_{(r)}G$ and
$C(M,\delta^{-1})\rtimes_{(r)}G$.
\end{thm}

\begin{ex}\label{ex-comp}
Assume that $L$ is a Lie group with finite component group such that
$L/K$, the quotient by the maximal compact subgroup $K$ of $L$, admits a
$L$-equivariant spin$^c$-structure (see the discussion at the beginning of
this section). Assume that $G$ is a cocompact discrete subgroup of $L$.
Then $G$ satisfies Axioms (A1), (A2), (A3) with respect to the
manifold $X=L/K$.

Let $P\subseteq L$
be any maximal parabolic subgroup.
Then $P$ is a closed connected
solvable subgroup of $L$ such that $M=L/P$
 is a compact manifold on which $L$, and hence $G$, acts
amenably (since $P$ is amenable, see \cite{AR}).

So the action of $G$ on $M$ satisfies all requirements of the above
theorem if the action of $G$ on $L/K \times L/P$ admits an equivariant
spin$^c$ structure. This will certainly be true if $L/K\times L/P$ admits
an $L$-invariant spin$^c$-structure, which easily follows if $L/P$
has an equivariant spin$^c$-structure for the action of $K$.
\end{ex}

\section{Extension of Poincar\'e Duality to Non-$G$-Compact Manifolds}\label{sec-nonGcompact}

In this section we want to show how to extend our main result  to
the case of non-$G$-compact manifolds $X$. Recall that for a
Hausdorff topological space $X$, the {\em $K$-homology of $X$ with
compact supports} is defined as  $K_*^c(X)=\lim_{Z} K_*(Z)$, where
$Z$ runs through the compact subsets of $X$ directed by inclusion.
Here we want to consider a non-commutative variant of this. If $G$
is a compact group, $X$ a locally compact $G$-space, $p:E\to X$ a
$G$-equivariant C*-algebra bundle over $X$, and $B$ a $G$-algebra,
we write
$$K^*_{c,G}(\Gamma_0(E); B):=\lim_Z \KK_*^G(\Gamma_0(E|_Z),B),$$
where $Z$ runs through the $G$-invariant compact subsets of $X$ and $E|_Z$ denotes the restriction of $E$ to $Z$.
We call $K^*_{c,G}(\Gamma_0(E); B)$
 the {\em $G$-equivariant $K$-homology of $\Gamma_0(E)$ with compact support in $X$
and coefficient $B$} and we simply write $K^*_{c,G}(\Gamma_0(E))$ if
$B=\CC$. In this section we shall prove

\begin{thm}\label{thm-noncompact} Let $G$ be a compact group which acts isometrically
on the complete Riemannian manifold $X$. Let $\delta\in \Br_G(X)$ be a smooth
element (the smoothness assumption can be avoided if $G$ is the trivial group).
Then, for any $G$-algebra $B$, there is a natural isomorphism
$$K^*_{c,G}(C_0(X,\delta^{-1}); B)\cong K_*^G(C_0(X,\delta\cdot\tau)\hot B).$$
In particular, for $B=\CC$ we obtain an isomorphism
$$K^*_{c,G}(C_0(X,\delta^{-1}))\cong K_*^G(C_0(X,\delta\cdot\tau))$$.
\end{thm}

\begin{remark}\label{rem-sym} We should point out that the above result is symmetric
in the sense that one can switch $\delta$ and $\delta^{-1}$ in the
formula (which is trivial), but one can also move the
Clifford-bundle to the other side, so that the first isomorphism
becomes $K^*_{c,G}(C_0(X,\delta^{-1}\cdot\tau); B)\cong
K_*^G(C_0(X,\delta)\hot B)$ (or $K^*_{c,G}(C_0(X,\delta\cdot\tau);
B)\cong K_*^G(C_0(X,\delta^{-1})\hot B).$)
\end{remark}

Another result of Poincar\'e duality which involves $K$-homology with compact supports
is given for proper actions of discrete groups.
In that case, if $p:E\to X$ is a C*-algebra bundle
over $X$ and $B$ is any C*-algebra (with trivial $G$-action), then we put
$$K^*_c(\Gamma_0(E)\rtimes G; B):=\lim_{Z}\KK_*(\Gamma_0(E|_Z)\rtimes G, B),$$
where $Z$ runs through the $G$-compact subsets of $X$. This is the
$K$-homology of $\Gamma_0(E)\rtimes G$ with compact support in $G\backslash X$
and coefficient $B$ (for the trivial group), if we consider $\Gamma_0(E)\rtimes G$ as
the section algebra of a C*-algebra bundle over $G\backslash X$ (which one can do
by \cite{Will}).
 We then extend Theorem \ref{thm-poincare-bundle} as follows:

\begin{thm}\label{thm-poincare-bundle-2} Suppose that $G$ is a discrete group acting properly on the
(complete) $G$-manifold $X$ and assume that $\delta\in \Br_G(X)$ is smooth.
 Then for every C*-algebra $B$ we get a natural isomorphism
$$K^*_c(C_0(X,\delta^{-1})\rtimes G; B) \cong K_*(C_0(X,\delta\cdot\tau)\rtimes G\hot B).$$
in particular, for $B=\CC$ we get $K^*_c(C_0(X,\delta^{-1})\rtimes
G)\cong K_*(C_0(X,\delta\cdot\tau)\rtimes G)$.
\end{thm}

Again, we should point out that one can move the symbol $\tau$ for
the Clifford bundle to the other side. The idea of the proof is to
use an exhaustive increasing sequence of open sub-manifolds $X_n$ of
$X$ with $G$-compact closures $\bar{X}_n$ such that $\bar{X}_n$ is a
manifold with boundary and then taking suitable limits over $n\in
\NN$. It is well-known to the experts that such a sequence exists,
but we shall give a short argument below. We are grateful to J\"org
Sch\"urmann, for providing the details and references for this
argument:

\begin{prop}\label{prop-schuer}
Let $G$ be a locally compact group acting properly by diffeomorphisms on the
second countable manifold $X$. Then the following are true:
\begin{enumerate}
\item $X$ admits a $G$-equivariant Riemannian metric (that is, $G$ acts
isometrically on $X$ with respect to this metric).
\item There exists an increasing sequence of open subspaces $(X_n)_{n\in\NN}$
of $X$ such that $X=\cup_{n\in\NN}X_n$ and such that $\bar{X}_n$ is a $G$-compact
manifold with boundary for all $n\in \NN$.
\end{enumerate}
\end{prop}
\begin{proof}
The existence of a $G$-equivariant metric is shown in
\cite[Theorem 4.2.4]{Pflaum}. This is all we have to get if $X$ is $G$-compact.
So assume from now on that $X$ is not $G$-compact.
In \cite[Theorem 4.2.4]{Pflaum} it is also shown that for every $G$-invariant
open cover of $X$, there exists a partition of unity consisting of $G$-invariant
differentiable maps subordinated to this cover. Since $G$ acts properly on $X$,
the quotient $G\backslash X$ is  also a second countable locally compact Hausdorff space,
and hence we can find an increasing sequence of open sets $U_n$ in $G\backslash X$
with compact closures $\bar{U}_n$
which covers $G\backslash X$. Taking inverse images $W_n:=q^{-1}(U_n)\subseteq X$,
we can then  find a partition of unity $(f_n)_{n\in \NN}$ subordinate to $(W_n)_{n\in\NN}$
with the above specified properties. Then
$$f:X\to [0,\infty); f(x)=\sum_{n\in \NN} n\cdot f_n(x)$$
is a $G$-invariant differentiable map, which is {\em $G$-proper},
which means that inverses of compact sets in $\RR$ are $G$-compact
sets in $X$. By Sard's theorem, we find an increasing sequence
$(r_n)_{n\in \NN}$ of regular values for $f$ such that $r_n\to
\infty$ for $n\to \infty$. Since $f$ is $G$-invariant, we can put
$\bar{X}_n=f^{-1}([0,r_n])$ and the sequence $(\bar{X}_n)_{n\in
\NN}$ then satisfies all requirements of (ii).
\end{proof}

\begin{defn}\label{def-decomposable}
A Riemannian manifold $(\bar{X},g)$ with boundary
$\partial X$ and with isometric $G$-action by some locally compact group $G$,
is said to be {\em decomposable near the boundary}
if there
exists  a new Riemannian
metric $g'$ on $\bar{X}$ satisfying:
\begin{itemize}
\item $G$ acts isometrically on
$(\bar{X},g')$.
\item There exists a neighbourhood of $\partial X$ in $(\bar{X},g')$
which is isometric to $(0,1]\times \partial X$ equipped with the
product metric of $g_{\R}$ and $g|_{\partial X}$.
\item The action of $G$ on $(0,1]\times \partial X$ is the product action,
where $G$ acts  trivially on $(0,1]$.
\end{itemize}
\end{defn}

 \begin{remark}\label{rem-Cliff}
Given two $G$-equivariant metrics $g$ and $g'$ on  $\bar{X}$, the Euclidean bundle
$(TX, g)$ is $G$-equivariantly isomorphic to $(TX, g')$, which follows from  the fact that
the two inner products in $T_xX$ differ by a completely positive transformation
which depends smoothly on $x\in X$, and hence induces a bundle isomorphism.
It follows directly from this that the Clifford bundle does not depend, up to $G$-isomorphism,
on the given choice of a $G$-equivariant metric on $\bar{X}$.
\end{remark}

Using standard methods from Differential Geometry, we can prove:

\begin{prop}\label{decomposition}
Let $(\bar{X},g)$ be a Riemannian manifold with boundary
$\partial X$ and with isometric $G$-action such  that $G\backslash \bar{X}$ is compact.
Then $(\bar{X},g)$ is  decomposable near the boundary.
 \end{prop}
 \begin{proof}
 Define a local diffeomorphism $\phi:
(-\infty,0]\times \partial X \longrightarrow \bar{X};\,\,(x,t)\mapsto
\Exp_xt\nu(x)$, where $\nu(x)$ is
the normal vector at $x$ to $\partial X$.
Since $t\mapsto \Exp_xt\nu(x)$  is the unique path from $x$ in
the direction of $\nu(x)$, and since the given metric $g$
on $\bar{X}$ is $G$-equivariant, this map is equivariant
with respect to the product action on $(-\infty,0]\times \partial X$
coming from the trivial action on $(-\infty,0]$ and the given action
on $ \partial X$.

For each $x\in \partial X$ there exists an open neighbourhood $W_x$ of
$x$ in $\partial X$ and $\eps_x>0$
such that  $\phi|_{(-\epsilon_x,0]\times W_x}$ is a diffeomorphism.
Since $\partial X$ is $G$-compact,
there exist finitely many $W_{x_1},\ldots,W_{x_n}$ such that
$GW_{x_1}\cup\cdots\cup GW_{x_n}=\partial X$.

Let $K=\partial X\times\partial X\setminus \cup_{h\in G}\cup_{i=1}^n
(hW_{x_i}\times hW_{x_i})$. Then for every $(x,y)\in K$,
$\rho(x,y)>0$ and, since $K$ is $G$-compact, $\min_{(x,y)\in
K}\rho(x,y)>0$. Let
$$\epsilon=\frac1{2}\min\{\min_{K}\rho(x,y),\,\epsilon_{x_1},\ldots, \eps_{x_n})\}>0.$$
Then $\phi: (-\eps,0]\times \partial X\to\bar{X}$ is a diffeomorphism onto its image.

One can now construct a
Riemannian metric $g'$ on $\bar{X}$ with isometric $G$-action
by a smooth convex combination
 such that the
restriction of $g'$ to $\bar{X}\setminus \phi((-\epsilon,0]\times\partial X)$ is
$g$ and the restriction of $g'$ to
$\phi([-\frac{\epsilon}{2},0] \times\partial X )$ is the image of
the product metric on $ [-\frac{\epsilon}{2},0] \times\partial X$.
Of course, after rescaling, we may replace $ [-\frac{\epsilon}{2},0]$ by $[0,1]$.

\end{proof}

If a Riemannian $G$-manifold $\bar{X}$ is decomposable near the boundary, we shall
assume from now on that the metric is given as in Definition \ref{def-decomposable}.
We then extend the notion of decomposability near the boundary to C*-algebra bundles
over $\bar{X}$ as follows:

\begin{defn}\label{def-decombundles}
Let $\bar{X}$ be a Riemannian $G$-manifold with boundary which is
decomposable near the boundary. A locally trivial C*-algebra bundle
$p:\bar{E}\to \bar{X}$ over $\bar{X}$
 is said to be {\em decomposable near the boundary}, if
 the restriction of $\bar{E}$ to $[0,1]\times \partial X$ (which is the closure of $(0,1]\times \partial X$ in $\bar{X}$)
 is $G$-equivariantly isomorphic
 to $[0,1]\times E_{0}$, where $E_{0}$ denotes the restriction of $\bar{E}$ to
 $\{0\} \times\partial X$.
 \end{defn}

 As in \S \ref{sec-twistorb} for the question of feasibility, it follows from
the general homotopy invariance of fibre bundles (see \cite{Huse})
that every locally
 trivial C*-algebra bundle $\bar{E}$ over $\bar{X}$ is decomposable near the boundary if $G$ is the trivial group.
 Moreover, as in Proposition \ref{lem-trans} for feasible bundles, we may conclude  that
 $p:\bar{E}\to\bar{X}$ is decomposable near the boundary whenever it is smooth and admits a $G$-equivariant Hermitian connection. In particular, this holds whenever
 \begin{itemize}
 \item $\bar{E}$ is smooth and
\item  $G$ is compact, or $G$ is discrete and acts properly on $\bar{X}$.
\end{itemize}
Using this notation we get

\begin{prop}\label{prop-boundary}
Suppose that $p:\bar{E}\to\bar{X}$ is a bundle
over the complete Riemannian $G$-manifold with boundary $\bar{X}$ such that
 $\bar{X}$ and $\bar{E}$  are decomposable near the boundary. Let $X$ denote the interior
of $\bar{X}$ and let $E$ denote the restriction of $\bar{E}$ to $X$.
Then
$$\iota^*: \RKK^G(\bar{X}; \Gamma_0(\bar{E})\otimes A, C_0(\bar{X}, B))
\to \RKK^G({X}; \Gamma_0({E})\otimes A, C_0({X}, B))$$ is an
isomorphism for all $G$-algebras $A$ and $B$, where
$\iota:X\to\bar{X}$ denotes inclusion.
\end{prop}

For the proof we need

\begin{lem}\label{lem}
Suppose that $X$ and $Y$ are locally compact $G$-spaces, $p:E\to X$ is a locally trivial C*-algebra bundle over $X$ with action of $G$ and $F:Y\times [0,1]\to X$ is a continuous map such that
$F^*E$ is $G$-equivariantly isomorphic to $F_0^*E\times[0,1]$ as bundles over $Y\times[0,1]$,
where $F_t:Y\to X$ denotes
the evaluation of $F$ at $t$ for all $t\in [0,1]$. Let
$$\varphi: \Gamma_0(F_0^*E)\to \Gamma_0(F_1^*E)$$
denote the isomorphism induced from the isomorphism $F_1^*E\cong
F_0^*E$ coming from the evaluation of the isomorphism $F^*E\cong
F_0^*E\times[0,1]$ at $Y\times\{1\}$. Then
$$F_0^*(\alpha)=\varphi^*(F_1^*(\alpha))\in \RKK^G(Y; \Gamma_0(F_0^*E)\hot A, C_0(Y, B))$$
for all $\alpha\in \RKK^G(X; \Gamma_0(E)\hot A, C_0(X,B))$.
\end{lem}
\begin{proof} Consider the element $F^*(\alpha)\in
\RKK^G(Y\times[0,1]; \Gamma_0(F^*E)\hot A, C_0(Y\times [0,1], B))$.
Forget the $C[0,1]$-structure to view $F^*(\alpha)$ as an element of
$\RKK^G(Y; \Gamma_0(F^*E)\hot A, C_0(Y\times [0,1], B))$. The
isomorphism $F^*E\cong F_0^*E\times[0,1]$ induces an isomorphism
$\Psi:\Gamma_0(F^*E)\stackrel{\cong}{\to} \Gamma_0(F_0^*E)\hot
C[0,1]$. We then get a $G$-equivariant and $C_0(Y)$-linear
$*$-homomorphism
$$\psi:\Gamma_0(F_0^*E)\to
\Gamma_0(F^*E);\quad
\psi(\xi)=\Psi^{-1}(\xi\otimes 1)$$
such that $\varphi$ is equal to the composition of $\psi$
 with evaluation $\ev_1:\Gamma_0(F^*E)\to \Gamma_0(F_1^*E)$.
Hence, the element
$$\psi^*(F^*(\alpha))\in \RKK^G(Y; \Gamma_0(F_0^*E)\hot A, C_0(Y\times [0,1], B)).$$
gives a homotopy in $\RKK^G(Y; \Gamma_0(F_0^*E)\hot A, C_0(Y, B))$
between  $F_0^*(\alpha)$ and $\varphi^*(F_1^*(\alpha))$.
\end{proof}

\begin{proof}[Proof of Proposition \ref{prop-boundary}]
Let
$\bar{X}_{\frac{1}{2}}:=\bar{X}\setminus((\frac{1}{2},1]\times\partial X)$. We then consider the
composition
$$\begin{CD}
 \RKK^G(\bar{X}; \Gamma_0(\bar{E})\otimes A, C_0(\bar{X}, B))\\
@V\iota_1^*VV\\
 \RKK^G({X}; \Gamma_0({E})\otimes A, C_0({X}, B))\\
 @V\iota_2^*VV\\
  \RKK^G(\bar{X}_{\frac{1}{2}}; \Gamma_0(\bar{E}_{\frac{1}{2}})\otimes A, C_0(\bar{X}_{\frac{1}{2}}, B))\\
@V\iota_3^*VV\\
 \RKK^G({X}_{\frac{1}{2}}; \Gamma_0({E}_{\frac{1}{2}})\otimes A, C_0({X}_{\frac{1}{2}}, B)),\\
 \end{CD}
 $$
where $X_{\frac{1}{2}}=\tilde{X}\setminus([\frac{1}{2},1]\times\partial X)$ denotes the interior of $\bar{X}_{\frac{1}{2}}$,
and
$\bar{E}_{\frac{1}{2}}$ and $E_{\frac{1}{2}}$ denote the restrictions of $\bar{E}$ to
$\bar{X}_{\frac{1}{2}}$ and ${X}_{\frac{1}{2}}$, respectively. We further let
$\iota=\iota_1:X\to\bar{X}$, $\iota_2:\bar{X}_{\frac{1}{2}}\to X$, and $\iota_3: X_{\frac{1}{2}}\to
\bar{X}_{\frac{1}{2}}$ denote the inclusions. Let $\bar{h}:\bar{X}\to \bar{X}_{\frac{1}{2}}$ denote the
homeomorphism defined by
$$\bar{h}(y)=y\quad\text{for $y\in \tilde{X}$}\quad\text{and}\quad \bar{h}(t,x)=
(\frac{t}{2},x)\quad\text{for $(t,x)\in [0,1]\times\partial X$}.$$
It is clear that $\bar{h}$ restricts to a homeomorphism $h:X\to X_{\frac{1}{2}}$.
The identity on $\bar{X}$ is homotopic to the
composition $\bar{X}\stackrel{\bar{h}}{\to}\bar{X}_{\frac{1}{2}}\stackrel{\iota_1\circ \iota_2}{\to}\bar{X}$
 by
$\bar{F}:\bar{X}\times [0,1]\to \bar{X}$ such that
\begin{align*}
\bar{F}(y,s)&=y\quad\quad\quad\quad\quad\text{for all $(y,s)\in \tilde{X}\times [0,1]$}\quad\text{and}\\
\bar{F}(t,x,s)&=(t-\frac{ts}{2},x)\quad\;\;\text{for $(t,x,s)\in[0,1]\times\partial X\times[0,1]$.}
\end{align*}
The restriction of $\bar{F}$ to $X\times[0,1]$ gives a homotopy $F:X\times[0,1]\to X$
between the identity on $X$ and the composition
$X\stackrel{h}{\to}X_{\frac{1}{2}}\stackrel{\iota_3\circ\iota_2}{\to} X$.
Since the bundle $p:\bar{E}\to \bar{X}$ is decomposable near the boundary $\partial X$, it is now
easy to see that
$$\bar{F}^*\bar{E}\cong \bar{F}_0^*\bar{E}\times[0,1]=\bar{E}\times [0,1]$$
and similarly that $F^*E\cong E\times[0,1]$. Indeed, since $\bar{F}_s$ is the identity on
$\tilde{X}:=\bar{X}\setminus((0,1]\times \partial X)$,
this structure is clear on the restriction of the bundle to $\tilde{X}\times[0,1]$, and on
the remaining part $[0,1]\times \partial X$ it follows from the decomposability
 of $\bar{E}$ near $\partial X$
that the restriction of $\bar{F}^*\bar{E}$ to $[0,1]\times \partial X\times [0,1]$ is isomorphic
to $[0,1]\times E_0\times [0,1]$, which is the restriction of $\bar{E}\times[0,1]$ to
$[0,1]\times \partial X\times [0,1]$.

Using this we can now apply the above lemma to see that $
\bar{h}^*\circ \iota_2^*\circ\iota_1^*$ is the identity on
$\RKK^G(\bar{X}; \Gamma_0(\bar{E})\otimes A, C_0(\bar{X}, B))$,
which proves that $\iota^*=\iota_1^*$ is injective. In a similar way
we see that $ h^*\circ \iota_3^*\circ\iota_2^*$ is the identity on $
\RKK^G({X}; \Gamma_0({E})\otimes A, C_0({X}, B))$.
Moreover, \mbox{$h^*:\RKK^G({X}_{\frac{1}{2}}; \Gamma_0({E}_{\frac{1}{2}})\otimes A,
C_0({X}_{\frac{1}{2}}, B))\to \RKK^G({X}; \Gamma_0({E})\otimes A,
C_0({X}, B))$} is an isomorphism with inverse $(h^{-1})^*$, since
$h:X\to X_{\frac{1}{2}}$ is a homeomorphism. A similar statement
holds for $\bar{h}^*$. Since $ h^*\circ \iota_3^*\circ\iota_2^*$ is
the identity we see that $\iota_3^*$ is surjective, and since
$\iota=\iota_1= \bar{h}^{-1}\circ \iota_3\circ h$ it follows from
this that $\iota^*$ is surjective, too.
\end{proof}

If $\bar{X}$ is a complete Riemannian manifold with interior $X$ which is decomposable near the
 boundary $\partial X$,
then it follows from the product structure
$ [0,1]\times \partial X$ near the boundary that the interior $X$ of $\bar{X}$
can be made into a complete manifold without boundary by changing the metric
on $[0,1)$ into a complete metric (i.e., via a suitable diffeomorphism $[0,1)\cong [0,\infty)$).
Since the action of $G$ on the interval is trivial, this change of metric is compatible with the
$G$-action.

We say that a locally trivial C*-algebra bundle $p:\bar{E}\to \bar{X}$
is  {\em feasible}, if its restriction to $X$ is feasible  in the sense of \S \ref{sec-twistorb}, where we
regard $X$ with the complete metric as discussed above.
Using this notation and  the above results together with
Theorem \ref{thm-kas-bundle} we get the following version for
manifolds with boundary

\begin{thm}\label{thm-kas-bundle1} Let $G$ be a locally compact group and
let $\bar{X}$ be a complete Riemannian $G$-manifold with boundary
 $\partial X$ which is decomposable near the boundary. Let $X$ denote the interior of $\bar{X}$
 and
suppose that $p:\bar{E}\to\bar{X}$ is a feasible
 locally trivial C*-algebra bundle which is decomposable near $\partial X$. Then, for
 each pair of $G$-algebras $A$ and $B$, there is a natural isomorphism
 $$\RKK^G(\bar{X};\Gamma_0(\bar{E})\otimes A, C_0(\bar{X}, B))\cong
 \KK^G(\Gamma_0(E)\hot_{C_0(X)}C_{\tau}(X)\hot A, B).$$
 \end{thm}

\begin{remark}\label{rem-triv} Recall from Proposition
\ref{decomposition} that $\bar{X}$ is always decomposable near the
boundary if $\bar{X}$ is $G$-compact. Recall also that by the
results of \S1 and the Appendix,
 the conditions on the bundle $p:\bar{E}\to\bar{X}$ are
always satisfied if $G$ is the trivial group and $\bar{E}$ is any
locally trivial C*-algebra bundle over $\bar{X}$, or if $G$ is
compact or discrete and acts properly on $\bar{X}$ and  $\bar{E}$
is a smooth $G$-bundle over $\bar{X}$ (which means in particular
that the action of $G$ on $\bar{E}$ is also smooth).  In particular,
if the compact or discrete group $G$ acts properly on $\bar{X}$, the
result applies for all smooth elements  $\delta\in \Br_G(\bar{X})$.
as introduced in \S  \ref{sec-prop}.
\end{remark}

\begin{proof}[Proof of Theorem \ref{thm-kas-bundle1}]
By Proposition \ref{prop-boundary} we get a natural isomorphism
$$\iota^*: \RKK^G(\bar{X}; \Gamma_0(\bar{E})\otimes A, C_0(\bar{X}, B))
\to \RKK^G({X}; \Gamma_0({E})\otimes A, C_0({X}, B))$$ and after
changing the metric on $X$ to make it complete, we can apply Theorem
\ref{thm-kas-bundle} to obtain a natural isomorphism
 $$\RKK^G({X};\Gamma_0({E})\otimes A, C_0({X}, B))\cong
 \KK^G(\Gamma_0(E)\hot_{C_0(X)}C_{\tau}(X)\hot A, B).$$
 This finishes the proof.
  \end{proof}

%


Suppose now that $X$ is a complete Riemannian manifold equipped with
a {\em proper} isometric action of the locally compact group $G$.
By Proposition \ref{prop-schuer} we can find an increasing sequence $(\bar{X}_n)_{n\in \NN}$
consisting of $G$-compact (hence complete) sub-manifolds with
boundary such that $X=\cup_{n\in \NN}X_n$, where $X_n$ denotes the
interior of $\bar{X}_n$, and by Proposition \ref{decomposition} we know
all $\bar{X}_n$ are decomposable near the boundary.
Therefore we may apply Theorem \ref{thm-kas-bundle1}
to obtain isomorphisms
\begin{equation}\label{thm-kas-bundle2}
\RKK^G(\bar{X}_n;\Gamma_0(\bar{E})\otimes A, C_0(\bar{X}_n, B))\cong
 \KK^G(\Gamma_0(E)\hot_{C_0(X_n)}C_{\tau}(X_n)\hot A, B).
\end{equation}
Let
$\delta$ be a smooth element in $\Br_G(X)$. For all $n\in \NN$ we
 denote by
$C_0(\bar{X}_n,\delta)$ the restriction of
$C_0(X,\delta)$ to $\bar{X}_n$ .

We are now ready to give

\begin{proof}[Proof of Theorem \ref{thm-noncompact}]
Let $(\bar{X}_n)_{n\in \NN}$ be as above.
Apply Theorem \ref{thm-kas-bundle1} to
$\Gamma_0(\bar{E})=C_0(\bar{X}_n,\delta)$ for each $n$ to get
$$\RKK^G(\bar{X}_n;C_0(\bar{X}_n,\delta)\otimes A, C_0(\bar{X}_n,
B))\cong
 \KK^G(C_0(X_n,\delta\cdot\tau)\hot A, B).$$
By Lemma \ref{lem-inverse} and the fact that $\bar{X}_n$ is compact,
we obtain
$$
\KK^G(A, C_0(\bar{X}_n,\delta^{-1})\hot
B)\stackrel{\cong}{\longrightarrow}
\KK^G(C_0(X_n,\delta\cdot\tau)\hot A, B).
$$
Using symmetry of Poincar\'e duality (Remark \ref{symmetry}), one
also has
$$
\KK^G(C_0(\bar{X}_n,\delta^{-1})\hot A,
B)\stackrel{\cong}{\longrightarrow} \KK^G(A,
C_0(X_n,\delta\cdot\tau)\hot B).
$$
Putting $A=\C$ and passing to the direct limit over $n$, we get
\begin{align*}
&K^*_{c,G}(C_0(X,\delta^{-1})\rtimes G;\,B)
=\lim_{n}\KK^G(C_0(\bar{X}_n,\delta^{-1}),
B)\\
&=\lim_n \KK^G(\C, C_0(X_n,\delta\cdot\tau)\hot B)
=K_*^G(C_0(X,\delta\cdot\tau)\hot B).
\end{align*}
\end{proof}

In a similar way we get

\begin{proof}[Proof of Theorem \ref{thm-poincare-bundle-2}]
Applying the maps $H$, $F$, $E$ in the proof of Theorem
\ref{thm-poincare-bundle} and the isomorphism from Theorem
\ref{thm-kas-bundle1} applied to
$\Gamma_0(\bar{E})=C_0(\bar{X}_n,\delta)$, we have the following
isomorphism for every pair of algebras $A$ and $B$ with trivial
$G$-action:
$$\KK(C_0(X_n,\delta\cdot\tau)\rtimes G\hot A,\,B)\cong \KK(A,\,C_0(\bar{X}_n,\delta^{-1})\rtimes G\hot
B),$$ where $X_n=\text{Int}(\bar{X}_n)$. By the symmetry of
Poincar\'e duality (Remark \ref{symmetry}), we can rewrite this into
$$\KK(C_0(\bar{X}_n,\delta^{-1})\rtimes G\hot A,\,B)\cong \KK(A,\,C_0(X_n,\delta\cdot\tau)\rtimes G\hot
B).$$ Putting $A=\C$ and passing to the direct limit over $n$, this completes the proof.
\end{proof}

\begin{remark}\label{rem-A} Notice that in the above proofs it is crucial to restrict to the case
$A=\CC$ before taking the limits over $n$, since $\KK$-theory is not continuous under taking direct limits.
\end{remark}

\section{Appendix}\label{append}

\begin{defn}\label{def-con} Let $p:E\to X$ be a
locally trivial, smooth bundle of C*-algebras over a manifold $X$.
Let $\Gamma^{\infty}(E)$ denote the algebra of smooth sections. By a
(Hermitian) \emph{connection} on $E$ we will understand a map
$$\nabla\colon \Gamma^\infty (TX)\hot \Gamma^{\infty}(E) \to \Gamma^{\infty}(E)$$ satisfying
\begin{enumerate}
\item $\nabla_{fV}(a) = f\nabla_V (a)$,
\item $\nabla_V(fa) = f\nabla_V (a)  + V(f)\nabla_V (a)$,
\item $\nabla_V (ab) = a\nabla_V(b) + \nabla_V (a) b.$
\item $\nabla_V(a^*) = \nabla_V(a)^*$,
\end{enumerate}
for all  $f\in C^{\infty}(X), V\in \Gamma^{\infty}(TX)$ and $a,b \in
\Gamma^{\infty}(E)$.
\end{defn}

Hermitian connections on \emph{trivial} bundles $E = X\times A$
always exist; for $p \in X$, $a \in C^\infty(X, A)$, define
$\nabla_V(a) := \frac{d}{dt}_{|_{t=0}} a\bigl(\gamma (t)\bigr)$ for
$\gamma\colon (-\epsilon , \epsilon ) \to X$ a smooth integral curve
for $V$. It is easy to check that the space of connections for $E$
is convex.  By a standard argument \cite{KN}
 it follows that a
Hermitian connection exists for any locally trivial bundle $E$.
Explicitly:
$$\nabla_V (s) = \sum_i  \rho_i \nabla^i_V(s),$$
where $\{\rho_i\}_{i\in I}$ is a suitable partition of unity for
$X$.

If $X$ is $G$-manifold and $E$ is a smooth $G$-equivariant bundle
(which requires also that $G$ acts via smooth automorphisms of the
bundle), then $G$ acts on the space of connections on $E$ by
$g\colon \nabla \mapsto \nabla^g$ with
$$\nabla^g_V(s) = g\bigl( \nabla_{g^{-1}V}(g^{-1}s)\bigr).$$
If $G$ is \emph{compact} then we may average an arbitrary Hermitian
connection to obtain a $G$-invariant one.

If $G$ is discrete and acts properly on $X$,
 then we can
write $X$ as the union of $G$-spaces $U_i$  each $G$-isomorphic to
$W_i \times_{H_i}G$, where $H_i$ is a finite subgroup of $G$ and
$W_i$ is an $H_i$-space. More precisely: $W_i$ is an open subset of
$X$, $H_i$ leaves $W_i$ invariant, and the $G$-saturation $G\cdot
W_i$ is a disjoint union of open subsets parameterized by the cosets
of $G/H_i$.

To obtain a $G$-invariant connection on $E$ it therefore suffices by
taking a $G$-invariant partition of unity subordinate to the cover
$\{U_i\}$, to find a $G$-invariant Hermitian connection on
$E|_{U_i}$ for each $U_i$. We may assume that $W_i$ is chosen
sufficiently small so that the restriction of $E$ to $W_i$ is
trivializable; we may then fix a connection on the restriction of
$E$ to $W_i$. Furthermore, by averaging over the finite group $H_i$,
we can assume that this connection is $H_i$-invariant. Denote it
$\overline{\nabla}$. Fix now a component $W_i'$ of $U_i$. There is
then $g \in G$ such that $g(W_i) = W_i'$, and hence $g$ maps the
restriction of $E$ to $W_i$ isomorphically to the restriction of $E$
to $W_i'$. We may use this isomorphism to construct a
$gH_ig^{-1}$-invariant connection on $E|_{W_i'}$.  Since
$\overline{\nabla}$ was $H_i$-invariant, the choice of $g$ is
immaterial.  In this way we obtain a $G$-invariant connection on
$E|_{U_i}$ as required.

We have proved:

\begin{lem}\label{lem-connection} Let $E$ be a locally trivial smooth
$G$-bundle of C*-algebras over a $G$-manifold $X$. Then, if either
$G$ is compact or discrete, there exists a $G$-invariant Hermitian
connection on $E$.
\end{lem}

Any two $G$-invariant connections on $E$ differ by a $G$-invariant
bundle map $TX \stackrel{\omega}{\longrightarrow} \mathrm{End}_\C(
\Gamma^{\infty}(E))$. Since $\omega (ab) = \omega (a)b + a\omega
(b)$,
 actually $\omega$ maps to the bundle $\mathrm{Der}(E)$
 of derivations of $E$ which are compatible with
 the adjoint operation: $\omega_V (a^*) = \omega_V(a)^*$. We
 emphasize that these derivations are \emph{bounded}, by the
 closed graph theorem.

Suppose that $E$ is as above and is equipped with a $G$-invariant
connection $\nabla$. Parallel sections $a$ along a smooth path in
$X$ are defined as usual with connections. This gives rise in the
usual fashion to the notion of parallel transport along a path. To
show that parallel transport exists, along a path, it suffices to
show this locally, since $X$ is connected. Since we always have the
trivial connection locally, we need to solve an ordinary
differential equation of the form $a'(t) = \omega_V(t) a(t).$ It is
obvious that a $G$-invariant connection gives rise to a
$G$-invariant parallel transport operation, and furthermore that
parallel transport respects the algebra multiplication: that is, if
$a$ and $b$ are parallel sections along a path, then so is $ab$,
because $\nabla$ satisfies the Leibnitz rule. So, putting everything
together, we get

\begin{prop}\label{lem-trans} Suppose that $X$ is a
$G$-manifold and $p:E\to X$ is a locally trivial smooth bundle of
C*-algebras over $X$. Suppose further that there exists a
$G$-invariant Hermitian connection for $E$. Then $E$ is feasible
in the sense of \S1.
In particular, if $G$ is compact or if
$G$ is discrete and acts properly on $X$, then every $G$-equivariant
smooth locally trivial bundle of C*-algebras over $X$ is feasible.
\end{prop}

\begin{proof}
Let $U$ and $P:U\times [0,1]\to X$ be as in Notations \ref{P}.
For $(x,y)\in U$ let $\gamma_{x,y}:[0,1]\to X$ be the unique
geodesic which joins $x$ and $y$ and let
$\varphi_{x,y,t}:E_{\gamma_{x,y}(t)}\to E_x$ denote the inverse of
the parallel transport $E_x\to E_{\gamma_{x,y}(t)}$ along the path
$\gamma_{x,y}$. Then $\varphi:P^*E\to P_0^*E$ defined fibre wise by
$\varphi_{x,y,t}$ is a $G$-equivariant bundle isomorphism.
\end{proof}

\def\mathcs{{\normalshape\text{C}}^{\displaystyle *}}

\end{document}